\documentclass{amsart}
\usepackage[utf8]{inputenc}
\usepackage{a4}
\usepackage{graphics}
\usepackage{color}
\usepackage{amssymb}
\usepackage[all]{xy}
\usepackage{amsmath}
\usepackage{stmaryrd}
\usepackage{subfigure}
\usepackage{longtable}
\usepackage{setspace}
\usepackage{booktabs}

\usepackage{pst-plot}
\usepackage{pstricks}
\usepackage{color}
\usepackage{nicefrac}

\CompileMatrices

\newtheorem{theorem}{Theorem}[section]

\newtheorem{proposition}[theorem]{Proposition}
\newtheorem{corollary}[theorem]{Corollary}

\theoremstyle{definition}

\newtheorem{definition}[theorem]{Definition}
\newtheorem{example}[theorem]{Example}

\newtheorem{construction}[theorem]{Construction}

\newtheorem{remark}[theorem]{Remark}
\theoremstyle{remark}

\numberwithin{equation}{section}
\def\Chi{{\mathbb X}}

\def\rel{{\rm rlv}}
\def\div{{\rm div}}

\def\quot{/\!\!/}
\def\mal{\! \cdot \!}

\def\reg{{\rm reg}}
\def\rq#1{\widehat{#1}}
\def\t#1{\widetilde{#1}}
\def\b#1{\overline{#1}}
\def\bangle#1{\langle #1 \rangle}

\def\KK{{\mathbb K}}
\def\TT{{\mathbb T}}
\def\ZZ{{\mathbb Z}}

\def\QQ{{\mathbb Q}}
\def\PP{{\mathbb P}}
\def\CDiv{\operatorname{CDiv}}
\def\WDiv{\operatorname{WDiv}}
\def\PDiv{\operatorname{PDiv}}

\def\Eff{{\rm Eff}}
\def\Mov{{\rm Mov}}

\def\Ample{{\rm Ample}}
\def\SAmple{{\rm SAmple}}

\def\ord{{\rm ord}}

\def\cov{{\rm cov}}

\def\Cl{\operatorname{Cl}}

\def\Pic{\operatorname{Pic}}

\def\Hom{{\rm Hom}}

\def\Supp{{\rm Supp}}
\def\Spec{{\rm Spec}}
\def\Proj{{\rm Proj}}

\def\cone{{\rm cone}}
\def\lin{{\rm lin}}

\def\topto#1{\stackrel{{\scriptscriptstyle #1}}{\longrightarrow}}

\def\cov{{\rm cov}}

\begin{document}
\title[Three lectures on Cox rings]%
{Three lectures on Cox rings}
\author[J.~Hausen]{J\"urgen Hausen} 
\address{Mathematisches Institut, Universit\"at T\"ubingen,
Auf der Morgenstelle 10, 72076 T\"ubingen, Germany}
\email{juergen.hausen@uni-tuebingen.de}


\maketitle

\section*{Introduction}

These are notes of an introductory course held 
at the conference 
``Torsors: Theory and Applications'' in Edinburgh, 
January 2011.
Cox rings play meanwhile an important role 
in Arithmetic and Algebraic Geometry, and 
in fact appeared independently in these 
fields~\cite{CTSaA,Co,HuKe}.
We present basic ideas and concepts, in 
particular, we treat the 
interaction with Geometric Invariant Theory.
The notes are kept as a survey; for 
details and proofs we refer to~\cite{ArDeHaLa}.

The first lecture begins with a rigorous 
definition of the Cox sheaf $\mathcal{R}$ 
of a variety~$X$ as a sheaf of algebras graded 
by the divisor class group $\Cl(X)$.
The Cox ring is then the algebra $\mathcal{R}(X)$ 
of global sections.
After recalling the basic correspondence between
graded algebras and quasitorus actions, we discuss
the relative spectrum $\rq{X} := \Spec_X \mathcal{R}$
of the Cox sheaf. We call $\rq{X} \to X$ the 
characteristic space; it coincides with the universal 
torsor~\cite{CTSaA,Sk} precisely for 
locally factorial varieties $X$.
The first lecture ends with a characterization of 
$\rq{X}$ in terms of Geometric Invariant Theory.

The aim of the second lecture is to present a 
machinery for encoding varieties via their Cox ring.
A first input is an explicit description of the 
variation of not necessarily quasiprojective good 
quotients for quasitorus actions on affine varieties.
This allows to encode torically embeddable varieties 
with finitely generated Cox ring in terms of what
we call ``bunched rings''.
Specialized to the case of toric varieties~\cite{Da,Oda,Fu}, 
the bunched ring data correspond to fans via Gale duality.
Moreover, the language of bunched rings extends basic 
features of techniques developed for weighted 
complete intersections~\cite{Do,IaFl} to the 
multigraded case. 
We show how to read off basic geometric properties 
from the defining data and briefly touch the relations
to Mori Theory.

The third Lecture is devoted to varieties with torus 
action.
Generalizing the case of a toric variety, we describe
the Cox ring in terms of the action.
In the case of a complexity one action, one obtains
a very explicit presentation in terms of trinomial 
relations.
Combining this with the language of bunched rings
leads to a concrete description of 
rational complete varieties with a complexity one action
turning them into an easy accessible class for 
concrete computations.
We demonstrate this in the case of $\KK^*$-surfaces.
Some classification results on Fano threefolds and 
del Pezzo surfaces based on this approach are included 
in the text.

\tableofcontents

I would like to thank the referee for his careful reading 
of the manuscript and many helpful comments.

\section{First lecture}

\subsection{Cox sheaves and Cox rings}
We work in the category of (reduced) varieties 
over an algebraically closed field
$\KK$ of characteristic zero.
The Cox ring $\mathcal{R}(X)$ will contain 
a lot of accessible information about 
the underlying variety~$X$,
the rough idea is to define it as
\begin{eqnarray*}
\mathcal{R}(X)
& := & 
\bigoplus_{[D] \in \Cl(X)} \Gamma(X,\mathcal{O}_X(D)),
\end{eqnarray*} 
where the grading is by the divisor class group 
$\Cl(X)$ of $X$.
A priori, it is not clear what the ring 
structure should be,
in particular if there is torsion in the 
divisor class group $\Cl(X)$.
Our first task is to clarify this; we closely
follow~\cite[Chap.~I]{ArDeHaLa}.

We begin with recalling the basic notions on 
divisors on the normal algebraic variety $X$. 
A {\em prime divisor\/} on $X$ is an irreducible 
hypersurface $D \subseteq X$.
The group of {\em Weil divisors\/} on $X$ is the 
free abelian group $\WDiv(X)$ generated by 
the prime divisors. 
We write $D \ge 0$ for a Weil divisor $D$, 
if it is a nonnegative linear combination 
of prime divisors. 
To every rational function 
$f \in \KK(X)^*$ one associates its 
{\em principal divisor\/}
$$
\div(f) 
\ := \
\sum_D \ord_D(f) D
\ \in \ 
\WDiv(X),
$$
where $D$ runs through the prime divisors of 
$X$ and $\ord_D(f)$ is the vanishing order 
of $f$ at $D$. The {\em divisor class group\/} 
$\Cl(X)$ of $X$ is the factor group of 
$\WDiv(X)$ by the subgroup $\PDiv(X)$ 
of principal divisors.
To every Weil divisor $D$ on $X$, 
one associates a {\em sheaf $\mathcal{O}_X(D)$\/} 
of $\mathcal{O}_X$-modules: for any
open $U \subseteq X$ one sets
\begin{eqnarray*}
\Gamma(U,\mathcal{O}_X(D))
& := & 
\{f \in \KK(X)^*; \; (\div(f) + D)_{\vert U} \ge 0 \}
\ \cup \
 \{0\},
\end{eqnarray*}
where the restriction map $\WDiv(X) \to \WDiv(U)$
is defined for a prime divisor $D$ as
$D_{\vert U} := D \cap U$ if it intersects $U$ 
and $D_{\vert U} := 0$ otherwise.
The sheaf $\mathcal{O}_X(D)$ is reflexive, 
i.e. canonically isomorphic to its double 
dual, and of rank one.
Note that for any two functions 
$f_1 \in \Gamma(U,\mathcal{O}_X(D_1))$ 
and 
$f_2 \in \Gamma(U,\mathcal{O}_X(D_2))$ 
the product $f_1f_2$ belongs to
$\Gamma(U,\mathcal{O}_X(D_1+D_2))$.

\begin{definition}
\index{divisorial algebra! sheaf of}
\index{sheaf of divisorial algebras}
The {\em sheaf of divisorial algebras\/} 
associated to a subgroup
$K \subseteq \WDiv(X)$
is the sheaf of $K$-graded 
$\mathcal{O}_X$-algebras
$$ 
\mathcal{S} 
\ := \ 
\bigoplus_{D \in K} \mathcal{S}_D,
\qquad\qquad
\mathcal{S}_D
\ := \ 
\mathcal{O}_X(D),
$$
where the multiplication in $\mathcal{S}$ 
is defined by multiplying homogeneous 
sections in the field of functions $\KK(X)$. 
\end{definition}

The sheaf of divisorial algebras associated 
to a finitely generated group of Weil divisors 
turns out to be a sheaf of normal algebras.
A crucial observation is the following.

\begin{proposition}
Let $\mathcal{S}$ be the sheaf of divisorial
algebras associated to a finitely generated 
group $K \subseteq \WDiv(X)$.
If the canonical map $K \to \Cl(X)$ sending 
$D \in K$ to its class $[D] \in \Cl(X)$ is 
surjective, then $\Gamma(X,\mathcal{S})$ is 
a unique factorization domain.
\end{proposition}

\begin{example}
On the projective line $X = \PP_1$,
consider $D := \{\infty\}$,
the group 
$K := \ZZ D$, 
and 
the associated $K$-graded sheaf 
of algebras $\mathcal{S}$.
Then we have isomorphisms
$$ 
\varphi_n 
\colon \KK[T_0,T_1]_n \ \to \ \Gamma(\PP_1, \mathcal{S}_{nD}),
\qquad
f \ \mapsto \ f(1,z),
$$
where $\KK[T_0,T_1]_n \subseteq \KK[T_0,T_1]$ 
denotes the vector space of all polynomials 
homogeneous of degree~$n$.
Putting them together we obtain a graded 
isomorphism
$$ 
\KK[T_0,T_1] 
\ \cong \
\Gamma(\PP_1, \mathcal{S}).
$$
\end{example}

\begin{example}[The affine quadric surface]
\label{ex:affquadsurf1}
Consider the two-dimensional affine variety
$$
X 
\ := \ 
V(\KK^3 ; \; T_1T_2 - T_3^2)
\ \subseteq \
\KK^3.
$$
We have the functions $f_i := T_{i \vert X}$
on $X$ and with the prime divisors 
$D_1 := V(X ; f_1)$ and $D_2 := V(X ; f_2)$
on $X$, we have 
$$ 
\div(f_1) \ = \ 2D_1,
\qquad
\div(f_2) \ = \ 2D_2,
\qquad
\div(f_3) \ = \ D_1+D_2.
$$
For $K := \ZZ D_1$, let $\mathcal{S}$ 
denote the associated sheaf of divisorial algebras.
Consider the sections 
$$
g_1 := 1 \ \in \ \Gamma(X,\mathcal{S}_{D_1}),
\qquad 
g_2 := f_3f_1^{-1} \ \in \ \Gamma(X,\mathcal{S}_{D_1}),
$$
$$
g_3 := f_1^{-1} \ \in \ \Gamma(X,\mathcal{S}_{2D_1}),
\qquad
g_4 := f_1 \ \in \ \Gamma(X,\mathcal{S}_{-2D_1}).
$$
Then $g_1,g_2$ generate $\Gamma(X,\mathcal{S}_{D_1})$
as a $\Gamma(X,\mathcal{S}_{0})$-module, 
and $g_3,g_4$ are inverse to each other. 
Moreover, we have 
$$ 
f_1 \ = \ g_1^2g_4,
\qquad
f_2 \ = \ g_2^2g_4,
\qquad
f_3 \ = \ g_1g_2g_4.
$$
Thus, $g_1,g_2,g_3$ and $g_4$ generate the 
$\KK$-algebra
$\Gamma(X,\mathcal{S})$. 
Setting $\deg(Z_i) := \deg(g_i)$, we obtain a 
$K$-graded isomorphism 
$$
\KK[Z_1,Z_2,Z_3^{\pm 1}] 
\ \to \ 
\Gamma(X,\mathcal{S}),
\qquad
Z_1 \mapsto g_1, \
Z_2 \mapsto g_2, \
Z_3 \mapsto g_3.
$$
\end{example}

We are ready to define the Cox ring.
From now on $X$ is a normal variety 
with $\Gamma(X,\mathcal{O}^*) = \KK^*$
and finitely generated divisor class 
group $\Cl(X)$. 
The idea is to start with the sheaf 
$\mathcal{S}$ of divisorial algebras 
associated to a group $K$ of Weil 
divisors projecting onto $\Cl(X)$ and
then to identify systematically isomorphic
homogeneous components $\mathcal{S}_D$ 
and $\mathcal{S}_D'$ via dividing by 
a suitable sheaf of ideals.

\begin{construction}
\label{constr:crtorsion}
\index{Cox ring!}\index{Cox sheaf!}%
Fix a subgroup 
$K \subseteq \WDiv(X)$ 
such that the map 
$c \colon K \to \Cl(X)$
sending $D \in K$ to its class 
$[D] \in \Cl(X)$
is surjective.
Let $K^0 \subseteq K$ be the kernel 
of $c$,
and let $\chi \colon K^0 \to \KK(X)^*$
be a character, i.e.~a group homomorphism, 
with
$$
\div(\chi(E))
\ = \ 
E,
\qquad
\text{for all } 
E \in K^0.
$$
Let $\mathcal{S}$ be the sheaf of divisorial algebras 
associated to $K$ and denote by $\mathcal{I}$ 
the sheaf of ideals of $\mathcal{S}$
locally generated by the sections
$1 - \chi(E)$, where $1$ is homogeneous
of degree zero, $E$ runs through $K^0$ and 
$\chi(E)$ is homogeneous of degree $-E$.
The {\em Cox sheaf\/} associated to 
$K$ and $\chi$ is the quotient sheaf 
$\mathcal{R} := \mathcal{S}/\mathcal{I}$
together with the $\Cl(X)$-grading
$$ 
\mathcal{R}
\ = \
\bigoplus_{[D] \in \Cl(X)}  \mathcal{R}_{[D]},
\qquad\qquad
\mathcal{R}_{[D]} 
\ := \ 
\pi \left( \bigoplus_{D' \in c^{-1}([D])} \mathcal{S}_{D'} \right).
$$
where $\pi \colon \mathcal{S} \to \mathcal{R}$ 
denotes the projection.
The Cox sheaf $\mathcal{R}$ is a quasicoherent sheaf 
of $\Cl(X)$-graded $\mathcal{O}_X$-algebras.
The {\em Cox ring\/} is the ring of global 
sections
$$ 
\mathcal{R}(X)
\ := \
\bigoplus_{[D] \in \Cl(X)}  \mathcal{R}_{[D]}(X),
\qquad\qquad
\mathcal{R}_{[D]}(X)
\ := \ 
\Gamma(X,\mathcal{R}_{[D]}).
$$
\end{construction}

In general, the Cox sheaf is not a sheaf 
of divisorial algebras.
However, the following shows that it is 
not too far from them.

\begin{remark}
\label{rem:crpops1}
Situation as in Construction~\ref{constr:crtorsion}.
Then, for every  $D \in K$, we have an 
isomorphism of sheaves
$$
\pi_{\vert \mathcal{S}_D} \colon \mathcal{S}_D \ \to \ \mathcal{R}_{[D]}.
$$
Moreover, for every open subset $U \subseteq X$, we have a 
canonical isomorphism on the level of sections
\begin{eqnarray*}
\Gamma(U,\mathcal{S}) 
 / 
\Gamma(U,\mathcal{I})
& \cong & 
\Gamma(U,\mathcal{S}/\mathcal{I}).
\end{eqnarray*}
In particular, the Cox ring $\mathcal{R}(X)$ is 
isomorphic to the quotient $\Gamma(X,\mathcal{S})/\Gamma(X,\mathcal{I})$ 
of global sections.
\end{remark}

\begin{proposition}
\label{prop:coxsheafunique}
If $K,\chi$ and $K',\chi'$ are data
as in Construction~\ref{constr:crtorsion},
then the associated Cox sheaves are isomorphic 
as graded sheaves.
\end{proposition}

\begin{example}
[The affine quadric surface, continued]
\label{ex:affquadsurf2}
Look again at the two-dimensional affine 
variety discussed in~\ref{ex:affquadsurf1}:
$$
X 
\ := \ 
V(\KK^3 ; \; T_1T_2 - T_3^2)
\ \subseteq \
\KK^3.
$$
The divisor class group $\Cl(X)$ is of order two; 
it is generated by $[D_1]$.
The kernel of the projection 
$K \to \Cl(X)$ is $K^0 = 2\ZZ D_1$ 
and a character as in 
Construction~\ref{constr:crtorsion} 
is 
$$ 
\chi \colon K^0 \ \to \ \KK(X)^*,
\qquad
2nD_1 \ \mapsto \ f_1^{n}.
$$
The ideal $\mathcal{I}$ is globally 
generated by the section $1 - f_1$, where 
$f_1 \in \Gamma(X,\mathcal{S}_{-2D_1})$.
Consequently, the Cox ring of $X$ 
is given as
$$
\mathcal{R}(X) 
\ \cong \ 
\Gamma(X,\mathcal{S})
/ 
\Gamma(X,\mathcal{I})
\ \cong \ 
\KK[Z_1,Z_2,Z_3^{\pm 1}]/ \bangle{1-Z_3^{-1}}
\ \cong \
\KK[Z_1,Z_2],
$$ 
where the $\Cl(X)$-grading on the polynomial 
ring $\KK[Z_1,Z_2]$ is given by 
$\deg(Z_1) = \deg(Z_2) = [D_1]$.
\end{example}

The construction of Cox sheaves 
(and thus also Cox rings) of a variety
can be made canonical by fixing a suitable 
point.
We say that $x \in X$ is {\em factorial}, 
if the local ring $\mathcal{O}_{X,x}$ is 
factorial; that means that every Weil divisor 
is principal near $x$.

\begin{construction}
\index{Cox ring! canonical}
\index{Cox sheaf! canonical}%
Let $x \in X$ be a factorial point
and consider the sheaf of divisorial algebras 
$\mathcal{S}^{x}$ associated to the 
group of Weil divisors avoiding $x$:
$$
K^{x} 
\ := \
\{D \in \WDiv(X); \ x \not\in \Supp(D)\}
\ \subseteq \WDiv(X).
$$
Let $K^{x,0} \subseteq K^{x}$ be the subgroup
consisting of principal divisors.
Then, for every $E \in K^{x,0}$, 
there is a unique 
$f_E \in \Gamma(X,\mathcal{S}_{-E})$,
which is defined near $x$ and satisfies
$$ 
\div(f_E) \ = \ E,
\qquad\qquad
f_E(x) \ = \ 1.
$$ 
The map $\chi^x \colon K^x \to \KK(X)^*$ 
sending $E$ to $f_E$ is a character 
as in~\ref{constr:crtorsion}.
We call the resulting Cox sheaf $\mathcal{R}^x$
the {\em canonical Cox sheaf of the pointed
space $(X,x)$.}
\end{construction}

As noted in~\ref{rem:crpops1}, the homogeneous 
components of the Cox sheaf are isomorphic 
to reflexive rank one sheaves $\mathcal{O}_X(D)$.
The following associates to every 
homogeneous section of the Cox ring a divisor, 
not depending on the choices made 
in~\ref{constr:crtorsion}.

\begin{construction}
\index{$[D]$-divisor!}%
In the setting of Construction~\ref{constr:crtorsion},
let $D \in K$ and consider 
an element $0 \ne f \in \mathcal{R}_{[D]}(X)$.
By Remark~\ref{rem:crpops1}, 
there is a (unique) 
$\t{f} \in \Gamma(X, \mathcal{S}_D)$ 
with $\pi(\t{f}) = f$.
The {\em $[D]$-divisor\/} and 
the {\em $[D]$-localization\/}
of  $f$ are
$$
\div_{[D]}(f) 
\ := \
\div(\t{f}) + D,
\qquad\qquad
X_{[D],f} \ := \ X \setminus \Supp(\div_{[D]}(f)).
$$ 
\end{construction}

We now discuss basic algebraic properties 
of the Cox ring.
Recall that a {\em Krull ring\/} is 
an integral ring $R$ with a family 
$(\nu_i)_{i \in I}$ of discrete valuations 
such that for every non-zero $f$ in the 
quotient field $Q(R)$, one has $\nu_i(f) \ne 0$ 
only for finitely many $i \in I$ 
and $f \in R$ if and only if $\nu_i(f) \ge 0$ 
holds for all $i \in I$.

\goodbreak

\begin{theorem}
\label{prop:crnormal}
The Cox ring $\mathcal{R}(X)$ is 
an (integral) normal Krull ring.
Moreover, one has the following 
statements on localization and units.
\begin{enumerate}
\item
For every non-zero homogeneous 
$f \in \mathcal{R}_{[D]}(X)$, 
there is a canonical isomorphism
\begin{eqnarray*}
\Gamma(X,\mathcal{R})_f
& \cong & 
\Gamma(X_{[D],f},\mathcal{R}).
\end{eqnarray*}
\item
Every homogeneous unit of $\mathcal{R}(X)$ is constant.
If $X$ is complete then even $\mathcal{R}(X)^* = \KK^*$
holds.
\end{enumerate}
\end{theorem}

Finally, we take a closer look at divisibility 
properties of Cox rings.
In general, they are not factorial, i.e.~unique 
factorization domains,
but the following weaker property will be guaranteed.

\begin{definition}
\label{def:factgrad}
\index{$K$-prime!}%
\index{$K$-height!}%
\index{factorially graded!}%
Consider an abelian group $K$ and 
a $K$-graded integral $\KK$-algebra
$R = \oplus_K R_w$.
\begin{enumerate}
\item
A non-zero non-unit $f \in R$ is {\em $K$-prime\/} if
it is homogeneous and $f | gh$ with homogeneous 
$g,h \in R$ implies $f | g$ or $f | h$.
\item
We say that $R$ is {\em factorially graded\/} 
if every homogeneous  non-zero non-unit $f \in R$
is a product of $K$-primes.
\end{enumerate}
\end{definition}

\begin{remark}
The concepts ``factorially graded'' and ``factorial'' 
coincide if the grading group is torsion free.
As soon as we have torsion in the grading group,
they may differ.
For example, consider
\begin{eqnarray*}
R 
& := &
\KK[T_1,T_2,T_3] / \bangle{T_1^2+T_2^2+T_3^2}.
\end{eqnarray*}
This is surely not a factorial ring.
But it becomes factorially 
graded by the group 
$K = \ZZ \oplus \ZZ/2\ZZ \oplus \ZZ/2\ZZ$ 
when we set
$$ 
\deg(T_1) \ := \ (1,\b{0},\b{0}),
\quad
\deg(T_2) \ := \ (1,\b{1},\b{0}),
\quad
\deg(T_3) \ := \ (1,\b{0},\b{1}).
$$
\end{remark}

\begin{theorem}
\label{thm:crfactgrad}
Suppose that the Cox ring $\mathcal{R}(X)$ 
satisfies $\mathcal{R}(X)^* = \KK^*$; 
for example, assume $X$ to be complete.
Then $\mathcal{R}(X)$ 
is factorially $\Cl(X)$-graded.
If moreover $\Cl(X)$ is torsion free, then 
$\mathcal{R}(X)$ is factorial.
\end{theorem}

The following statement shows that divisibility 
in the Cox ring $\mathcal{R}(X)$ can be formulated
geometrically in terms of $[D]$-divisors on the underlying 
variety~$X$.

\begin{proposition}
Suppose that the Cox ring $\mathcal{R}(X)$ 
satisfies $\mathcal{R}(X)^* = \KK^*$.
\begin{enumerate}
\item
An element $0 \ne f \in \Gamma(X,\mathcal{R}_{[D]})$ 
divides  $0 \ne g \in \Gamma(X,\mathcal{R}_{[E]})$ 
if and only if $\div_{[D]}(f) \le \div_{[E]}(g)$
holds.
\item 
An element $0 \ne f \in  \Gamma(X,\mathcal{R}_{[D]})$ 
is $\Cl(X)$-prime
if and only if the divisor $\div_{[D]}(f) \in \WDiv(X)$ 
is prime. 
\end{enumerate}
\end{proposition}


\subsection{Quasitorus actions}
Here we recall the correspondence between 
affine algebras $A$ graded by a finitely 
generated abelian group $K$ 
and affine varieties $X$ coming with 
an action of a quasitorus $H$.
Moreover, we discuss good quotients and their 
basic properties.
Standard references are~\cite{Bo,OnVi,Sp,Kr,MuFoKi}.

An {\em affine algebraic group\/} 
is an affine variety $G$ together with a group 
structure such that the group operations are 
morphisms. A {\em morphism of affine algebraic groups\/}
is a morphism of the underlying varieties which
is moreover a group homomorphism.
A {\em character\/} is a morphism $G \to \KK^*$
to the multiplicative group of the ground field.
Endowed with pointwise multiplication, the 
characters of $G$ form an abelian group $\Chi(G)$.

\begin{definition}
A {\em quasitorus}, also called a 
{\em diagonalizable group\/}, 
is an affine algebraic group $H$ 
whose algebra of regular functions 
$\Gamma(H,\mathcal{O})$ is generated 
as a $\KK$-vector space by the characters
$\chi \in \Chi(H)$.
A {\em torus\/} is a connected quasitorus.
\end{definition}

\begin{example}
The {\em standard $n$-torus\/} is
$\TT^n := (\KK^*)^n$.
Its characters are precisely the Laurent monomials
$T^{\nu} = T_1^{\nu_1} \cdots T_n^{\nu_n}$, 
where $\nu \in \ZZ^n$, and its 
algebra of regular functions is the Laurent polynomial 
algebra
$$ 
\Gamma(\TT^n,\mathcal{O})
\ = \ 
\KK[T_1^{\pm 1}, \ldots,T_n^{\pm 1}]
\ = \ 
\bigoplus_{\nu \in \ZZ^n} \KK \cdot T^{\nu}
\ = \ 
\KK[\ZZ^n].
$$
\end{example}

To any finitely generated abelian group $K$ 
one associates in a functorial way a
quasitorus, namely $H := \Spec\, \KK[K]$, 
the spectrum of the group algebra $\KK[K]$.

\begin{remark}
The quasitorus $H:= \Spec\, \KK[K]$ 
can be realized as a closed 
subgroup of a standard $r$-torus as follows.
By the elementary divisors theorem,
we find generators $w_1, \ldots, w_r$ of $K$ 
such that the epimorphism $\ZZ^r \to K$, 
$e_i \mapsto w_i$ has kernel 
$$
\ZZ a_1e_1 \oplus \ldots \oplus \ZZ a_se_{s}
\ \subseteq \
\ZZ^r,
\qquad 
a_1, \ldots, a_s \in \ZZ_{\ge 1}.
$$
The corresponding morphism $H \to \TT^r$
is a closed embedding realizing $H$ as the 
kernel of $\TT^r \to \TT^s$,
$(t_1, \ldots, t_r) \mapsto (t_1^{a_1}, \ldots, t_s^{a_s})$.
In particular, we see that $H$ is a direct 
product of a torus and a finite abelian 
group:
$$
H 
\ \cong \
C(a_1)  \times  \ldots \times C(a_s) \times \TT^{r-s},
\qquad\qquad
C(a_i) 
\ := \ 
\{\zeta \in \KK^*; \; \zeta^{a_i} = 1\}.
$$
\end{remark}

\begin{theorem}
We have contravariant exact functors
being essentially inverse to each other:
\begin{eqnarray*}
\{\text{finitely generated abelian groups}\}
& \longleftrightarrow &
\{\text{quasitori}\}
\\
K 
& \mapsto & 
\Spec \, \KK[K],
\\
\Chi(H)
& \mapsfrom &
H.
\end{eqnarray*}
Under these equivalences, 
the free finitely generated abelian groups 
correspond to the tori.
\end{theorem}

Quasitori are linearly reductive in characteristic 
zero: every rational representation even splits into 
one-dimensional ones. Applying this to the 
representation on the algebra of global functions
of an affine variety with quasitorus action, we 
obtain a grading by the character group in the 
following way.

\begin{remark}
Let a quasitorus $H$ act on a
not necessarily affine variety~$X$.
Then the algebra $\Gamma(X, \mathcal{O})$ becomes 
$\Chi(H)$-graded via
$$
\Gamma(X, \mathcal{O})
\ = \ 
\bigoplus_{\chi \in \Chi(H)} \Gamma(X, \mathcal{O})_\chi, 
\qquad 
\Gamma(X, \mathcal{O})_\chi
\ := \ 
\{f \in \Gamma(X, \mathcal{O}); \;  f(h \, \mal \, x) = \chi(h) f(x) \}. 
$$
\end{remark}

We now associate in functorial manner 
to every affine algebra $A = \oplus_K A_w$ 
graded by a finitely generated abelian group 
$K$ the affine variety $X = \Spec \, A$
with an action of the quasitorus 
$H = \Spec \, \KK[K]$.
Again this can be made concrete.

\begin{construction}
Let $K$ be a finitely generated abelian group
and $A$ a $K$-graded affine algebra. 
Set $X = \Spec \, A$.
If $f_i \in A_{w_i}$, $i=1,\ldots,r$, generate $A$, 
then we have a closed embedding
$$ 
X \ \to \ \KK^r,
\qquad
x \ \mapsto \ (f_1(x), \ldots, f_r(x)),
$$
and $X \subseteq \KK^r$ is invariant under 
the diagonal action of $H = \Spec \, \KK[K]$
given by the characters $\chi^{w_1}, \ldots, \chi^{w_r}$.
Note that for any $f \in A$ homogeneity is 
characterized by 
\begin{eqnarray*}
f \in A_w
& \iff & 
f (h \mal x) = \chi^w(h) f(x)
\text{ for all } h \in H, \, x \in X.
\end{eqnarray*}
\end{construction}

\begin{theorem}
\label{prop:gradalg2afftorac}
We have contravariant functors
being essentially inverse to each other:
\begin{eqnarray*}
\{\text{graded affine algebras}\}
& \longleftrightarrow &
\{\text{affine varieties with quasitorus action}\}
\\
A
& \mapsto & 
\Spec \, A,
\\
\Gamma(X,\mathcal{O})
& \mapsfrom &
X.
\end{eqnarray*}
\end{theorem}

\begin{definition}
\index{good quotient!}%
\index{quotient! good}%
\index{geometric quotient!}%
\index{quotient! geometric}%
\label{def:goodquot}
Let $G$ be a linearly reductive affine 
algebraic group $G$ acting on a variety $X$.
A morphism $p \colon X \to Y$ is called 
a {\em good quotient\/} for this action  
if it has the following properties:
\begin{enumerate}
\item
$p \colon X \to Y$ is affine and $G$-invariant,
\item
the pullback $p^* \colon \mathcal{O}_Y \to (p_*\mathcal{O}_X)^G$
is an isomorphism.
\end{enumerate}
A morphism $p \colon X \to Y$ is called a 
{\em geometric quotient\/} if it is  a good quotient
and its fibers are precisely the $G$-orbits.
\end{definition}

\begin{remark}
Let a linearly reductive group $G$ act on
an affine variety $X$.
Then Hilbert's finiteness theorem ensures that 
the algebra of invariants 
$\Gamma(X,\mathcal{O})^G$ is finitely 
generated.
This gives a good quotient $p \colon X \to Y$
with $Y := \Spec \, \Gamma(X,\mathcal{O})^G$.
\end{remark}

\begin{remark}
Let $A = \oplus_K A_w$ be an affine algebra graded 
by a finitely generated abelian group $K$.
Then $H = \Spec \, \KK[K]$ acts on $X = \Spec \, A$
and we have 
\begin{eqnarray*}
\Gamma(X,\mathcal{O})^H
& = & 
A_0.
\end{eqnarray*}
In order to compute $A_0$, choose homogeneous generators 
$f_1, \ldots, f_r$ of $A$ and consider the homomorphism
$$ 
Q \colon \ZZ^r \ \to \ K, 
\qquad 
e_i \ \mapsto \deg(f_i).
$$
Then, for any set $B$ of generators of the monoid
$\ZZ^r_{\ge 0} \cap  \ker(Q)$, the algebra $A_0$
of invariants is generated by the products 
$f_1^{\nu_{1}} \cdots f_r^{\nu_{r}}$, where $\nu \in B$.
\end{remark}

The basic properties of good quotients
are that they map closed invariant subsets 
to closed sets and that they separate 
disjoint closed invariant sets.
An immediate consequence is that the
target space carries the quotient 
topology.
Another application is the following 
statement on the fibers.

\begin{proposition}
\label{prop:goodquotfibers}
Let a linearly reductive algebraic group $G$ 
act on a variety $X$, and let 
$p \colon X \to Y$ be a good quotient.
Then $p$ is surjective and 
for any $y \in Y$ one has:
\begin{enumerate}
\item
There is exactly one closed $G$-orbit 
$G \mal x$ in the fiber $p^{-1}(y)$.
\item
Every orbit $G \mal x' \subseteq p^{-1}(y)$
has $G \mal x$ in its closure.
\end{enumerate}
\end{proposition}

The first statement means that a good 
quotient $p \colon X \to Y$ parametrizes 
the closed orbits of the $G$-variety $X$. 
Using the description of the fibers one 
easily verifies that a good quotient is
categorical, i.e.~universal with respect to 
invariant morphisms. 
In particular, the quotient space is unique
up to isomorphism which justifies the common
notation $X \quot G$.

\begin{example}
\label{exam:1torusonk21}
Consider the $\KK^*$-action $t \cdot (z_1,z_2) = (t^az_1,t^bz_2)$ 
on $\KK^2$.
The following three cases are typical.
\begin{enumerate}
\item
We have $a=b=1$. Every $\KK^*$-invariant function is 
constant and the constant map $p \colon \KK^2 \to \{{\rm pt}\}$
is a good quotient.
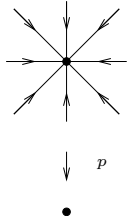
\begin{center}
\begin{picture}(0,0)%
\includegraphics{elliptic.pstex}%
\end{picture}%
\setlength{\unitlength}{829sp}%
\begingroup\makeatletter\ifx\SetFigFontNFSS\undefined%
\gdef\SetFigFontNFSS#1#2#3#4#5{%
  \reset@font\fontsize{#1}{#2pt}%
  \fontfamily{#3}\fontseries{#4}\fontshape{#5}%
  \selectfont}%
\fi\endgroup%
\begin{picture}(3666,6454)(1318,-6482)
\put(4051,-5011){\makebox(0,0)[lb]{\smash{{\SetFigFontNFSS{6}{7.2}{\familydefault}{\mddefault}{\updefault}{\color[rgb]{0,0,0}$p$}%
}}}}
\end{picture}%

\end{center}
\item
We have $a=0$ and $b=1$. 
The algebra of $\KK^*$-invariant functions
is generated by $z_1$ and  the map $p \colon \KK^2 \to \KK$,
$(z_1,z_2) \mapsto z_1$ is a good quotient.
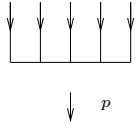
\begin{center}
\begin{picture}(0,0)%
\includegraphics{parabolic.pstex}%
\end{picture}%
\setlength{\unitlength}{829sp}%
\begingroup\makeatletter\ifx\SetFigFontNFSS\undefined%
\gdef\SetFigFontNFSS#1#2#3#4#5{%
  \reset@font\fontsize{#1}{#2pt}%
  \fontfamily{#3}\fontseries{#4}\fontshape{#5}%
  \selectfont}%
\fi\endgroup%
\begin{picture}(3846,4566)(1228,-4594)
\put(4051,-3211){\makebox(0,0)[lb]{\smash{{\SetFigFontNFSS{6}{7.2}{\familydefault}{\mddefault}{\updefault}{\color[rgb]{0,0,0}$p$}%
}}}}
\end{picture}%

\end{center}
\item
We have $a=1$ and $b=-1$. 
The algebra of $\KK^*$-invariant functions
is generated by $z_1z_2$ and  $p \colon \KK^2 \to \KK$,
$(z_1,z_2) \mapsto z_1z_2$ is a good quotient.
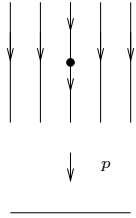
\begin{center}
\begin{picture}(0,0)%
\includegraphics{hyperbolic.pstex}%
\end{picture}%
\setlength{\unitlength}{829sp}%
\begingroup\makeatletter\ifx\SetFigFontNFSS\undefined%
\gdef\SetFigFontNFSS#1#2#3#4#5{%
  \reset@font\fontsize{#1}{#2pt}%
  \fontfamily{#3}\fontseries{#4}\fontshape{#5}%
  \selectfont}%
\fi\endgroup%
\begin{picture}(3846,6366)(1228,-6394)
\put(4051,-5011){\makebox(0,0)[lb]{\smash{{\SetFigFontNFSS{6}{7.2}{\familydefault}{\mddefault}{\updefault}{\color[rgb]{0,0,0}$p$}%
}}}}
\end{picture}%

\end{center}
Note that the general $p$-fiber is a single 
$\KK^*$-orbit, whereas $p^{-1}(0)$ consists
of three orbits and is reducible.
\end{enumerate}
\end{example}

\begin{example}[The affine quadric surface as a quotient] 
\label{ex:affquadsurf3}
Consider the action of the multiplicative group 
$C(2) = \{1,-1\}$ on $\KK^2$ given by  
\begin{eqnarray*}
\zeta \cdot z
& := & 
(\zeta z_1, \zeta z_2).
\end{eqnarray*}
This action has $0$ as a fixed point and any 
$z \ne 0$ has trivial isotropy group.
The algebra $A \subseteq \KK^2$ of invariants 
is generated by 
$$
f_{11} := T_1^2, \qquad
f_{22} := T_2^2, \qquad
f_{12} := T_1T_2.
$$
The ideal of relations among them is generated by 
the polynomial $f_{11}f_{22}-f_{12}^2$.
Consequently, with  
$X := V(w_1w_2-w_3^2) \subseteq \KK^3$,
the quotient map is
$$ 
\pi \colon \KK^3 \ \to \ X, 
\qquad
z \ \mapsto \ (f_{11}(z),f_{22}(z),f_{12}(z)).
$$
The fibers of $\pi$ are precisely the 
$C(2)$-orbits. In particular, $\pi$ is a
geometric quotient (as it holds 
for any finite group action on an affine 
variety).
\end{example}

\begin{example}[The affine quadric threefold as a quotient] 
\label{ex:affquadthreef1}
Consider the action of $\KK^*$ on $\KK^4$ given 
by 
\begin{eqnarray*}
t \cdot z
& := & 
(t z_1, t^{-1} z_2, t z_3, t^{-1}z_4).
\end{eqnarray*}
Note that $\KK^*$ has $0$ as a fixed point and
any $z \ne 0$ has trivial isotropy group.
The algebra $A \subseteq  \KK[T_1,\ldots,T_4]$ 
of invariants is generated by 
$$
f_{12} := T_1T_2, \qquad
f_{34} := T_3T_4,  \qquad
f_{14} := T_1T_4, \qquad
f_{23} := T_2T_3.
$$ 
The ideal of relations is generated by 
the polynomial $f_{12}f_{34} - f_{14}f_{23}$. 
Thus, with 
$X := V(w_1w_2-w_3w_4) \subseteq \KK^4$,
the quotient map is
$$ 
\pi \colon \KK^4 \ \to \ X, 
\qquad
z \ \mapsto \ (f_{12}(z),f_{34}(z),f_{14}(z),f_{23}(z)).
$$
The fiber over any point $0 \ne x$ is a free 
$\KK^*$-orbit and hence isomorphic to $\KK^*$, 
whereas the fiber over the point $0 \in X$ is
given by
$$ 
\pi^{-1}(0) 
\ = \ 
V(T_1,T_3) \cup V(T_2,T_4) 
\ \subseteq \ 
\KK^4.
$$
\end{example}


\subsection{Characteristic spaces}
As before, $X$ is a normal variety with 
$\Gamma(X,\mathcal{O}^*) = \KK^*$.
We discuss the geometric counterpart 
of a Cox sheaf $\mathcal{R}$ on $X$,
its relative spectrum.
In order to obtain a reasonable object, 
$\mathcal{R}$ should be locally of finite
type.
By Theorem~\ref{prop:crnormal}, this holds, 
if the Cox ring $\mathcal{R}(X)$ is finitely 
generated.
Moreover, one can show that $\mathcal{R}$ 
is locally of finite type for every 
$\QQ$-factorial $X$.

\begin{construction}
\label{constr:univtors}
\index{space!characteristic}%
\index{characteristic space!}%
\index{characteristic quasitorus!}%
\index{quasitorus! characteristic}%
Let $\mathcal{R}$ be a Cox sheaf on $X$
and suppose that $\mathcal{R}$ is locally 
of finite type.
Then the relative spectrum
\begin{eqnarray*}
\rq{X}
& := & 
\Spec_X(\mathcal{R})
\end{eqnarray*}
is a quasiaffine variety.
The $\Cl(X)$-grading of the sheaf $\mathcal{R}$ 
defines an action of the diagonalizable group
\begin{eqnarray*}
H_X 
& := & 
\Spec \, \KK[\Cl(X)]
\end{eqnarray*}
on $\rq{X}$. The canonical morphism 
$q_X \colon \rq{X} \to X$
is a good quotient for this action,
and we have an isomorphism of graded sheaves
\begin{eqnarray*}
\mathcal{R}
& \cong &
(q_X)_* (\mathcal{O}_{\rq{X}}). 
\end{eqnarray*}
We call $q_X \colon \rq{X} \to X$ the 
{\em characteristic space\/} associated 
to $\mathcal{R}$, 
and $H_X$ the
{\em characteristic quasitorus\/} of~$X$. 
\end{construction}

In the case of a locally factorial variety $X$, 
the characteristic space coincides with the 
universal torsor over $X$. 
As soon as $X$ has non-factorial singularities 
the two concepts differ from 
each other, as we will indicate below.

\begin{proposition}
\label{prop:etalequot}%
\label{prop:smallmap2}%
\label{cor:quasaf}%
Consider the  characteristic space 
$q_X \colon \rq{X} \to X$.
\begin{enumerate}
\item
The inverse image $q_X^{-1}(X_{\reg})$
of the set of smooth points 
is smooth, 
$H_X$ acts freely there
and $q_X^{-1}(X_{\reg}) \to X_{\reg}$
is an \'etale $H_X$-principal bundle. 
\item
For any closed set $A \subseteq X$ 
of codimension at least two, 
the inverse image 
$q_X^{-1}(A) \subseteq \rq{X}$
is as well of codimension at least two.
\item
Let $\rq{x} \in \rq{X}$ be a point such that 
the orbit $H_X \mal \rq{x} \subseteq \rq{X}$ is closed,
and consider an element $f \in \Gamma(X,\mathcal{R}_{[D]})$.
Then we have 
\begin{eqnarray*}
f(\rq{x}) \ = \ 0
& \iff &
q_X(\rq{x}) \ \in \ \Supp(\div_{[D]}(f)).
\end{eqnarray*}
\end{enumerate}
\end{proposition}

We now relate properties of the $H_X$-action 
to geometric properties on $X$.
For $x \in X$, let $\PDiv(X,x) \subseteq \WDiv(X)$ 
denote the subgroup of all Weil divisors, which 
are principal on some neighbourhood of~$x$. 
We define the {\em local class group\/} 
of $X$ at~$x$ to be the factor group
\begin{eqnarray*}
\Cl(X,x) 
& := & 
\WDiv(X) / \PDiv(X,x).
\end{eqnarray*}
Obviously the group $\PDiv(X)$ of principal 
divisors is contained in $\PDiv(X,x)$. 
Thus, there is a canonical epimorphism 
$\pi_x \colon \Cl(X) \to \Cl(X,x)$.
We denote by $H_{X,\rq{x}} \subseteq H_X$
the isotropy group of $\rq{x} \in \rq{X}$.

\begin{proposition}
\label{prop:isogrp2divisors2}
Consider the characteristic space 
$q_X \colon \rq{X} \to X$.
Given~$x \in X$, fix a point
$\rq{x} \in q_X^{-1}(x)$ with 
closed $H_X$-orbit. Then we have
a canonical isomorphism
\begin{eqnarray*}
\Cl(X,x) & \cong & \Chi(H_{X,\rq{x}}).
\end{eqnarray*}
\end{proposition}

Recall that a point $x \in X$ is factorial 
if and only if near $x$ every Weil divisor 
is principal.
We say that $x \in X$ is {\em $\QQ$-factorial\/} 
if near $x$ for every Weil divisor some multiple 
is principal.

\goodbreak

\begin{corollary}
Consider the characteristic space 
$q_X \colon \rq{X} \to X$.
\begin{enumerate}
\item
A point $x \in X$ is factorial if and only if 
the fiber $q_X^{-1}(x)$ is a single $H_X$-orbit with 
trivial isotropy.
\item
A point $x \in X$ is $\QQ$-factorial if and only if 
the fiber $q_X^{-1}(x)$ is a single $H_X$-orbit.
\end{enumerate}
\end{corollary}

The first part of the following says in particular
that characteristic space and universal torsor coincide 
if and only if $X$ has at most factorial singularities.

\begin{corollary}
\label{cor:univtorquot}
Consider the characteristic space 
$q_X \colon \rq{X} \to X$.
\begin{enumerate}
\item
The action of $H_X$ on $\rq{X}$ is free if and only if 
$X$ is locally factorial.
\item
The good quotient $q_X \colon \rq{X} \to X$ is geometric 
if and only if $X$ is $\QQ$-factorial.
\end{enumerate}
\end{corollary}

Recall that the {\em Picard group\/} 
of $X$ is the factor group of the group $\CDiv(X)$ 
of locally principal Weil divisors by the subgroup 
of principal divisors:
$$
\Pic(X) 
\ = \ 
\CDiv(X) / \PDiv(X)
\ = \ 
\bigcap_{x \in X} \ker(\pi_x).
$$

\begin{corollary}
\label{prop:picviacoxrings}
Consider the characteristic space 
$q_X \colon \rq{X} \to X$.
Let $\rq{H}_X \subseteq H_X$ 
be the subgroup generated by
all isotropy groups 
$H_{X,\rq{x}}$, where $\rq{x} \in \rq{X}$. 
Then we have 
\begin{eqnarray*}
\ker \bigl(\Chi(H_X) \to \Chi(\rq{H}_X)\bigr)
& = & 
\bigcap_{\rq{x} \in \rq{X}} \ker\bigl(\Chi(H_X) \to \Chi(H_{X,\rq{x}})\bigr)
\end{eqnarray*}
and the projection $H_X \to H_X/\rq{H}_X$ 
corresponds to the inclusion 
$\Pic(X) \subseteq \Cl(X)$ of character 
groups. 
\end{corollary}

\begin{corollary}
\label{cor:fixpt2pic0}
If the variety $\rq{X}$ contains an $H_X$-fixed point,
then the Picard group $\Pic(X)$ is trivial.
\end{corollary}

As we noted, the characteristic space of a 
variety~$X$ is a quasiaffine variety~$\rq{X}$ 
with an action of the characteristic 
quasitorus $H_X$ having $X$ as a good 
quotient. 
Our next aim is to characterize this situation
in terms of Geometric Invariant Theory.

\begin{definition}
\label{def:stronglystable}
\index{strongly stable action!}%
\index{action! strongly stable}%
Let $G$ be an affine algebraic group and 
$W$ a $G$-variety. 
We say that the $G$-action on $W$ is 
{\em strongly stable\/} if there is an 
open invariant subset $W' \subseteq W$ 
with the following properties:
\begin{enumerate}
\item
the complement $W \setminus W'$ is of codimension
at least two in $W$,
\item 
the group $G$ acts freely, i.e.~with 
trivial isotropy groups, on $W'$,
\item
for every $x \in W'$ the orbit $G \mal x$ 
is closed in $W$.
\end{enumerate}
\end{definition}

\begin{definition}
Let a group $G$ act on a normal 
variety $Y$.
A divisor $\sum a_DD$ is called 
{\em $G$-invariant\/}
if its multiplicities satisfy 
$a_{g \mal D} = a_D$ for every $g \in G$.
We say that $Y$ is {\em $G$-factorial\/} 
if every $G$-invariant divisor is principal.
\end{definition}

Note that a quasiaffine variety with a quasitorus 
$H$ acting on it is $H$-factorial if and only 
if its ring of functions is factorially graded 
by the character group $\Chi(H)$.

\begin{remark}
Let $X$ be a normal variety with 
characteristic space
$q_X \colon \rq{X} \to X$.
Then $\rq{X}$ is $H_X$-factorial
and $q_X^{-1}(X_{\reg}) \subseteq \rq{X}$
satisfies the required properties 
of $W' \subseteq W$ of Definition~\ref{def:stronglystable}.
\end{remark}

\begin{theorem}
\label{bigfreequot}
Let a quasitorus $H$ act on a normal quasiaffine 
variety $\mathcal{X}$ with a good quotient
$q \colon \mathcal{X} \to X$.
Assume that 
$\Gamma(\mathcal{X},\mathcal{O}^*) = \KK^*$ holds, 
$\mathcal{X}$ is $H$-factorial
and the $H$-action is strongly stable.
Then there is a commutative diagram
$$ 
\xymatrix{
{\mathcal{X}}
\ar[dr]_{q}
\ar[rr]^{\mu}_{\cong}
&&
{\rq{X}}
\ar[dl]^{q_{X}}
\\
&
X
&
}
$$
where the quotient space $X$ is a normal variety 
with $\Gamma(X,\mathcal{O}_X^*) = \KK^*$,
we have $\Cl(X) = \Chi(H)$ and 
$q_X \colon \rq{X} \to X$ is a characteristic space
for $X$ and the isomorphism 
$\mu \colon \mathcal{X} \to \rq{X}$ is
equivariant with respect to the actions of $H = H_X$.
\end{theorem}

\begin{example}[The affine quadric surface, continued]
In~\ref{ex:affquadsurf3}, we realized the surface  
$X = V(T_1T_2-T_3^2)$ as a strongly stable 
quotient of $\KK^2$ by $\ZZ/2\ZZ$.
Thus, $\Cl(X)$ is of order two and $\KK^2 \to X$ 
is a characteristic space.
Moreover, from~\ref{cor:fixpt2pic0} we infer $\Pic(X) = 0$.
\end{example}

\begin{example}[The affine quadric threefold, continued]
In~\ref{ex:affquadthreef1}, we realized the affine threefold
$X = V(T_1T_2-T_3T_4)$ as a strongly stable 
quotient of $\KK^4$ by $\KK^*$.
Thus, we have $\Cl(X) \cong \ZZ$ and $\KK^4 \to X$ 
is a characteristic space.
Moreover, from~\ref{cor:fixpt2pic0} we infer $\Pic(X) = 0$.
\end{example}

These two examples are special cases of
the much bigger class of  
{\em toric varieties\/}, i.e., normal 
varieties $X$ endowed with a torus action
$T \times X \to X$ such that for some 
$x_0 \in X$ the orbit map $T \to X$,
$t \mapsto t \mal x_0$ is an open embedding.
Toric varieties admit a complete description
in terms of lattice fans; standard references
are~\cite{Da,Oda,Fu, CoLiS}.
In this picture, Cox ring and characteristic 
space look as follows; see~\cite{Co,Au,Ba,Mus}.

\begin{construction}
\label{constr:toriccox}%
Assume that the toric variety $X$ arises from a 
fan $\Sigma$ in a lattice $N$.
The condition $\Gamma(X,\mathcal{O}^*) = \KK^*$ 
means that the primitive vectors 
$v_1, \ldots, v_r \in N$ on the rays of $\Sigma$ 
generate $N_\QQ$ as a vector space.
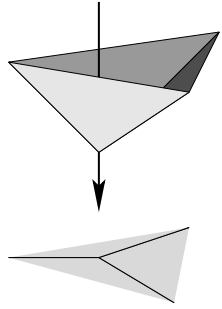
\begin{center}
\begin{picture}(0,0)%
\includegraphics{cox.pstex}%
\end{picture}%
\setlength{\unitlength}{829sp}%
\begingroup\makeatletter\ifx\SetFigFontNFSS\undefined%
\gdef\SetFigFontNFSS#1#2#3#4#5{%
  \reset@font\fontsize{#1}{#2pt}%
  \fontfamily{#3}\fontseries{#4}\fontshape{#5}%
  \selectfont}%
\fi\endgroup%
\begin{picture}(6366,9087)(868,-8194)
\end{picture}%

\end{center}
Set $F := \ZZ^r$ and consider 
the linear map $P \colon F \to N$ sending the 
$i$-th canonical base vector $f_i \in F$
to $v_i \in N$.
There is a fan $\rq{\Sigma}$ in $F$ 
consisting of certain faces of the positive 
orthant 
$\delta \subseteq F_{\QQ}$, namely
$$ 
\rq{\Sigma}
\ := \
\{\rq{\sigma} \preceq \delta; \; P(\rq{\sigma}) \subseteq \sigma
\text{ for some } \sigma \in \Sigma\}.
$$
The fan $\rq{\Sigma}$ defines an open toric 
subvariety $\rq{X}$ of  
$\b{X} = \Spec(\KK[\delta^\vee \cap E])$,
where $E := \Hom(F,\ZZ)$.
Note that all rays $\cone(f_1), \ldots, \cone(f_r)$
of the positive orthant $\delta \subseteq F_{\QQ}$
belong to $\rq{\Sigma}$ and thus we have 
$$ 
\Gamma(\rq{X},\mathcal{O})
\ = \ 
\Gamma(\b{X},\mathcal{O})
\ = \ 
\KK[\delta^\vee \cap E].
$$
As $P \colon F \to N$ is a map of the fans 
$\rq{\Sigma}$ and $\Sigma$, i.e., sends cones
of $\rq{\Sigma}$ into cones of $\Sigma$,
it defines a morphism $p \colon \rq{X} \to X$
of toric varieties. 
Note that we have $\b{X} = \KK^r$ and in terms 
of the coordinates $T_i = \chi^{e_i}$, the open subset 
$\rq{X} \subseteq \b{X}$ is given as 
$$
\rq{X} 
\ = \
\b{X} \setminus V(T^\sigma; \; \sigma \in \Sigma),
\qquad
T^\sigma \ = \ T_1^{\varepsilon_1} \cdots T_r^{\varepsilon_r},
\quad
\varepsilon_i = 
\begin{cases}
1, & v_i \not\in\sigma,
\\
0, & v_i \in\sigma.
\end{cases}
$$
Now, consider the dual map 
$P^* \colon M \to E$, where $M := \Hom(N,\ZZ)$,
set $K := E/P^*(M)$ and denote by 
$Q \colon E \to K$ the projection.
We obtain a $K$-grading of the polynomial 
ring $\KK[T_1, \ldots, T_r]$ by setting
$$
\deg(T_i) \ := \ Q(e_i) \ \in \ K.
$$
This gives an action of
$H := \Spec \, \KK[K]$ on $\b{X}$.
The set $\rq{X} \subseteq \b{X}$ is 
invariant and $p \colon \rq{X} \to X$ 
is a good quotient.
Moreover, 
$$
W \ := \ 
\b{X} \setminus \bigcup_{i \ne j} V(T_i,T_j)
\ \subseteq \ 
\rq{X}
$$
satisfies the conditions 
of~Definition~\ref{def:stronglystable} 
for this action. 
Thus, $p \colon \rq{X} \to X$ 
is a characteristic space for $X$. 
In particular, we have
$$ 
\Cl(X) \ \cong \ K,
\qquad
\mathcal{R}(X) \ \cong \ \KK[T_1, \ldots, T_r].
$$
\end{construction}


\section{Second lecture}

\subsection{Variation of good quotients}
Given a variety $X$ with an action 
of a linearly reductive group $G$,
the task of Geometric Invariant Theory 
is to describe the {\em good $G$-sets},
i.e.~the invariant open subsets 
$U \subseteq X$ admitting a good quotient 
$U \to U \quot G$. In general, there will
be several such good $G$-sets; this effect 
is also called ``variation of good quotients''.
For details, see~\cite[Sec.~III.1]{ArDeHaLa}
and the original references~\cite{MuFoKi,BeHa3,ArHa3}.

\begin{example}
\label{exam:1torusonk22}
For the $\KK^*$-action $t \cdot (z_1,z_2) = (t^az_1,t^bz_2)$ 
on $\KK^2$; as in~\ref{exam:1torusonk21}, we consider
the three typical cases.
\begin{enumerate}
\item
Let $a=b=1$. Besides the good quotient 
$\KK^2 \to \{ {\rm pt} \}$, 
there is a nice geometric quotient 
$\KK^2 \setminus \{0\} \to \PP_1$,
$(z_1,z_2) \mapsto [z_1,z_2]$.

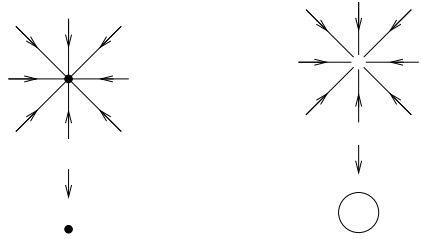
\begin{center}
\begin{picture}(0,0)%
\includegraphics{ell1.pstex}%
\end{picture}%
\setlength{\unitlength}{829sp}%
\begingroup\makeatletter\ifx\SetFigFontNFSS\undefined%
\gdef\SetFigFontNFSS#1#2#3#4#5{%
  \reset@font\fontsize{#1}{#2pt}%
  \fontfamily{#3}\fontseries{#4}\fontshape{#5}%
  \selectfont}%
\fi\endgroup%
\begin{picture}(3666,6454)(1318,-6482)
\end{picture}%

\qquad\qquad\qquad
\begin{picture}(0,0)%
\includegraphics{ell2.pstex}%
\end{picture}%
\setlength{\unitlength}{829sp}%
\begingroup\makeatletter\ifx\SetFigFontNFSS\undefined%
\gdef\SetFigFontNFSS#1#2#3#4#5{%
  \reset@font\fontsize{#1}{#2pt}%
  \fontfamily{#3}\fontseries{#4}\fontshape{#5}%
  \selectfont}%
\fi\endgroup%
\begin{picture}(3666,6955)(1318,-6983)
\end{picture}%

\end{center} 
\item
Let $a=0$ and $b=1$. Besides the  
good quotient $\KK^2 \to \KK$, $(z_1,z_2) \mapsto z_1$,
there is a geometric quotient
$\KK^2  \setminus V(T_2) \to \KK$,
$(z_1,z_2) \mapsto z_1$.

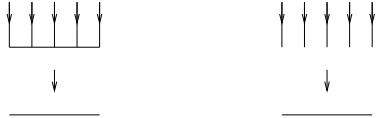
\begin{center}
\begin{picture}(0,0)%
\includegraphics{par1.pstex}%
\end{picture}%
\setlength{\unitlength}{622sp}%
\begingroup\makeatletter\ifx\SetFigFont\undefined%
\gdef\SetFigFont#1#2#3#4#5{%
  \reset@font\fontsize{#1}{#2pt}%
  \fontfamily{#3}\fontseries{#4}\fontshape{#5}%
  \selectfont}%
\fi\endgroup%
\begin{picture}(3846,4566)(1228,-4594)
\end{picture}

\qquad\qquad\qquad
\begin{picture}(0,0)%
\includegraphics{par2.pstex}%
\end{picture}%
\setlength{\unitlength}{622sp}%
\begingroup\makeatletter\ifx\SetFigFont\undefined%
\gdef\SetFigFont#1#2#3#4#5{%
  \reset@font\fontsize{#1}{#2pt}%
  \fontfamily{#3}\fontseries{#4}\fontshape{#5}%
  \selectfont}%
\fi\endgroup%
\begin{picture}(3846,4566)(1228,-4594)
\end{picture}

\end{center} 
\item
Let $a=1$ and $b=-1$. Besides the good quotient
$\KK^2 \to \KK$, $(z_1,z_2) \mapsto z_1z_2$, there are 
two geometric ones:
$\KK^2  \setminus V(T_i) \to \KK$,
$(z_1,z_2) \mapsto z_1z_2$.

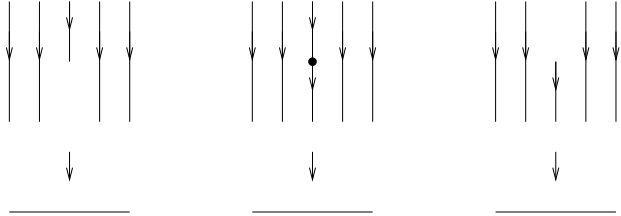
\begin{center}
\begin{picture}(0,0)%
\includegraphics{hyp0.pstex}%
\end{picture}%
\setlength{\unitlength}{829sp}%
\begingroup\makeatletter\ifx\SetFigFontNFSS\undefined%
\gdef\SetFigFontNFSS#1#2#3#4#5{%
  \reset@font\fontsize{#1}{#2pt}%
  \fontfamily{#3}\fontseries{#4}\fontshape{#5}%
  \selectfont}%
\fi\endgroup%
\begin{picture}(3846,6366)(1228,-6394)
\end{picture}%

\qquad\qquad
\begin{picture}(0,0)%
\includegraphics{hyp1.pstex}%
\end{picture}%
\setlength{\unitlength}{829sp}%
\begingroup\makeatletter\ifx\SetFigFontNFSS\undefined%
\gdef\SetFigFontNFSS#1#2#3#4#5{%
  \reset@font\fontsize{#1}{#2pt}%
  \fontfamily{#3}\fontseries{#4}\fontshape{#5}%
  \selectfont}%
\fi\endgroup%
\begin{picture}(3846,6366)(1228,-6394)
\end{picture}%

\qquad\qquad
\begin{picture}(0,0)%
\includegraphics{hyp2.pstex}%
\end{picture}%
\setlength{\unitlength}{829sp}%
\begingroup\makeatletter\ifx\SetFigFontNFSS\undefined%
\gdef\SetFigFontNFSS#1#2#3#4#5{%
  \reset@font\fontsize{#1}{#2pt}%
  \fontfamily{#3}\fontseries{#4}\fontshape{#5}%
  \selectfont}%
\fi\endgroup%
\begin{picture}(3846,6366)(1228,-6394)
\end{picture}%

\end{center} 
\end{enumerate}
\end{example}

All the good $G$-sets occurring in this 
example are maximal w.r.t. inclusions
of the following type:
a subset $U' \subseteq U$
of a good $G$-set $U \subseteq X$ is
called {\em $G$-saturated\/} in $U$
if it satisfies $U' = \pi^{-1}(\pi(U'))$,
where $\pi \colon U \to U \quot G$
is the good quotient.
Any good $G$-set is $G$-saturated 
in a maximal one and for a description
it is reasonable to focus on the latter ones.

Our aim is to present concrete combinatorial
descriptions for quasiprojective and for 
torically embeddable quotients of quasitorus 
actions on certain affine varieties.
Let us fix the setting.
By $K$ we denote a finitely generated 
abelian group and we consider 
an affine $K$-graded $\KK$-algebra
\begin{eqnarray*}
A 
& = & 
\bigoplus_{w \in K} A_w.
\end{eqnarray*}
\index{orbit monoid!}%
\index{orbit group!}%
Then the quasitorus $H := \Spec \, \KK[K]$
acts on the affine variety $X := \Spec \, A$.
Let $K_{\QQ} := K \otimes_{\ZZ} \QQ$ denote the 
rational vector space associated to $K$.
Given $w \in K$, we write again $w$ for the 
element $w \otimes 1 \in K_{\QQ}$.

\begin{definition}
The {\em weight cone\/} of the $K$-graded algebra 
$A$ is the convex polyhedral cone
$$
\omega_X 
\ := \ 
\omega(A)
\ = \ 
\cone(w \in K; \;  A_w \ne \{0\})
\ \subseteq \ 
K_{\QQ}.
$$
To every point $x \in X$, we associate its 
{\em orbit cone\/}, this is the convex polyhedral 
cone 
$$ 
\omega_x 
\ := \ 
\cone(w \in K; \; f(x) \ne 0 \text{ for some } f \in A_w) 
\ \subseteq \ 
\omega_X.
$$
\end{definition}

In order to see that these cones are indeed 
polyhedral, let $f_1, \ldots, f_r$ be homogeneous 
generators for $A$ and set $w_i := \deg(f_i)$. 
Then the weight cone $\omega_X$ is generated by 
$w_1, \ldots, w_r$ and the orbit cone $\omega_x$
is generated by those $w_i$ with $f_i(x) \ne 0$.
In particular, we see that $\omega_X$ is the 
general orbit cone and that there are only finitely 
many orbit cones.

\goodbreak

\begin{example}
\label{exam:1torusonk23}
We determine the orbit cones of $\KK^*$-action 
$t \cdot (z_1,z_2) = (t^az_1t^bz_2)$ on $\KK^2$;
again, we consider the three typical cases:
\begin{enumerate}
\item
We have $a=b=1$. 
The weight cone is $\omega_{\KK^2} = \QQ_{\ge 0}$ 
and the possible orbit cones are 
$$
\omega_{(0,0)} = \{0\},
\qquad
\omega_{(1,1)} = \QQ_{\ge 0}.
$$
\item
We have $a=0$ and $b=1$. 
The weight cone is $\omega_{\KK^2} = \QQ_{\ge 0}$ 
and the possible orbit cones are 
$$
\omega_{(0,0)} = \{0\}, 
\qquad
\omega_{(1,0)} = \QQ_{\ge 0}.
$$
\item
We have $a=1$ and $b=-1$. 
The weight cone is $\omega_{\KK^2} = \QQ$ 
and the possible orbit cones are 
$$
\omega_{(0,0)} = \{0\},
\quad 
\omega_{(1,0)} = \QQ_{\ge 0},
\quad
\omega_{(0,1)} = \QQ_{\le 0},
\quad
\omega_{(1,1)} = \QQ.
$$
\end{enumerate}
\end{example}

In this example we considered a subtorus action 
on a toric variety, and we used the fact that 
it suffices to determine one orbit cone for each 
toric orbit. This idea can be generalized
using equivariant embeddings and gives 
the following concrete recipes for computing
orbit cones.

\begin{remark}
\label{constr:ffgeo}
Fix a system of homogeneous generators 
$\mathfrak{F} = (f_{1}, \dots, f_{r})$ for 
our $K$-graded algebra $A$.
Then $H$ acts diagonally on $\KK^r$ via the 
characters $\chi^{w_1}, \ldots, \chi^{w_r}$,
where $w_i = \deg(f_i)$,
and we have an $H$-equivariant closed embedding
$$ 
X \ \to \ \KK^r,
\qquad
x \ \mapsto \ (f_1(x), \ldots, f_r(x)).
$$
With $E = \ZZ^r$ and $\gamma = \cone(e_1, \ldots, e_r)$, 
we may identify $\KK[T_1, \ldots, T_r]$ with 
$\KK[E \cap \gamma]$ and thus 
regard $\KK^r$ as the affine toric variety
associated to the dual cone $\delta := \gamma^\vee$.
For any $\gamma_0 \preceq \gamma$ and 
$\delta_0 := \gamma_0^\perp \cap \delta$, the following 
statements are equivalent:
\begin{enumerate}
\item
the product over all $f_{i}$ with $e_{i} \in \gamma_{0}$
lies not in 
$\sqrt{\bangle{f_{j}; \; e_{j} \not\in \gamma_{0}}} \subseteq A$,
\item
there is a point $z \in X$ with
$z_i \ne 0 \Leftrightarrow e_i \in \gamma_0$
for all $1 \le i \le r$,
\item 
the toric orbit $\TT^r \mal z_{\delta_0} \subseteq \KK^r$ 
corresponding to $\delta_0 \preceq \delta$ 
meets $X$,
\item
the intersection $\delta_0^\circ \cap {\rm Trop}(X)$ 
with the tropical variety is non-empty.
\end{enumerate}
In order to determine the orbit cones, denote by
$Q \colon E \to K$ the homomorphism sending 
$e_i$ to $w_i$. 
Moreover, we call {\em $\mathfrak{F}$-faces}
the $\gamma_0 \preceq \gamma$ satisfying~(i).
Then the orbit cones of $X$ are precisely 
the images $Q(\gamma_0)$, where $\gamma_0 \preceq \gamma$
is an $\mathfrak{F}$-face.
\end{remark}

\begin{example}
\label{ex:delpezzo1}
Set $K := \ZZ^2$ and consider the $K$-grading of
$\KK[T_1, \ldots, T_5]$ defined by $\deg(T_i) := w_i$,
where $w_i$ is the $i$-th column of the matrix
\begin{eqnarray*}
Q
& := &
\left[
\begin{array}{rrrrr}
1 & -1 & 0 & -1 & 1
\\
1 & 1 & 1 & 0 & 2
\end{array}
\right].
\end{eqnarray*}
The corresponding action of $H = \TT^2$ on $\KK^5$
leaves $X := V(T_1T_2 + T_3^2 + T_4T_5)$ invariant.
The possible orbit cones $\omega_x$ for $x \in X$
are 
$$ 
\{0\}, \quad
\cone(w_1), \quad
\cone(w_2), \quad
\cone(w_4), \quad
\cone(w_5), 
$$
$$
\cone(w_1,w_4), \quad
\cone(w_2,w_4), \quad
\cone(w_1,w_5), \quad
\cone(w_2,w_5). 
$$
\end{example}

\begin{definition}
\index{GIT-cone!}%
The {\em GIT-cone\/} 
of an element $w \in \omega_X$
is the (nonempty) intersection 
of all orbit cones containing it:
\begin{eqnarray*} 
\lambda(w) 
& := &
\bigcap_{{x \in X,} \atop {w \in \omega_x}} \omega_x.
\end{eqnarray*}
We write $\Lambda(X,H)$ for the set of all
GIT-cones.
The {\em set of semistable points\/} 
associated to a GIT-cone 
$\lambda \subseteq \omega_X$ is
$$
X^{ss}(\lambda)
\ = \ 
\{x \in X; \; \lambda \subseteq  \omega_x\}
\ \subseteq  \ 
X.
$$
\end{definition}

\begin{remark}
Given a GIT-cone $\lambda \in \Lambda(X,H)$ and any weight 
$w \in \lambda^\circ$ in its relative interior, one easily 
checks
\begin{eqnarray*}
X^{ss}(\lambda)
& = &
\{x \in X; \; f(x) \ne 0 \text{ for some } f \in A_{nw}, \, n > 0\}.
\end{eqnarray*}
That means that $X^{ss}(\lambda)$ is the set of 
semistable points associated to the linearization 
of the trivial bundle given by the character 
$\chi^w$ in the sense of Mumford.
\end{remark}

\begin{example}
\label{exam:1torusonk24}
We compute GIT-cones and associated sets of 
semistable points for the $\KK^*$-action 
$t \cdot (z,w) = (t^az,t^bw)$ on $\KK^2$ in the
three typical cases.
$$
\xymatrix@R=3pt{
a=1, b=1:
&
{\lambda(0) = \{0\}} 
&
{X^{ss}(\lambda(0)) = \KK^2},
\\
&
{\lambda(1) = \QQ_{\ge 0}}
&
{X^{ss}(\lambda(1)) = \KK^2 \setminus \{0\}},
\\
\\
a=0, b=1:
&
{\lambda(0) = \{0\}} 
&
{X^{ss}(\lambda(0)) = \KK^2},
\\
&
{\lambda(1) = \QQ_{\ge 0}}
&
{X^{ss}(\lambda(1)) = \KK^2 \setminus V(T_2)},
\\
\\
a=-1, b=1:
&
{\lambda(0) = \{0\}} 
&
{X^{ss}(\lambda(0)) = \KK^2},
\\
&
{\lambda(-1) = \QQ_{\le 0}}
&
{X^{ss}(\lambda(-1)) =  \KK^2 \setminus V(T_1)},
\\
&
{\lambda(1) = \QQ_{\ge 0}}
&
{X^{ss}(\lambda(1)) = \KK^2 \setminus V(T_2)}.
\\
\\
}
$$ 
\end{example}

In the following statement, we mean by a {\em quasifan\/} 
a finite collection $\Lambda$ 
of not necessarily pointed polyhedral cones in a rational vector 
space such that any two $\lambda_1,\lambda_2 \in \Lambda$ 
intersect in a common face and for $\lambda \in \Lambda$ 
also every face of $\lambda$ belongs to~$\Lambda$.

\begin{theorem}
\label{thm:gitfan}
The collection $\Lambda(X,H) = \{\lambda(w); \; w \in \omega_X\}$
of all GIT-cones is a quasifan in $K_{\QQ}$
having the weight cone $\omega_X$ as its support.
\begin{enumerate}
\item
For every $\lambda \in \Lambda(X,H)$, there is a good 
quotient $X^{ss}(\lambda) \to Y(\lambda)$ for the 
action of $H$ on $X^{ss}(\lambda)$.
\item
For any two GIT-cones $\lambda_1, \lambda_2 \in \Lambda(X,H)$, 
we have $\lambda_2 \preceq \lambda_1$ if and only if
$X^{ss}(\lambda_1) \subseteq X^{ss}(\lambda_2)$ holds.
\item
If $X^{ss}(\lambda_1) \subseteq X^{ss}(\lambda_2)$ holds,
then there is an induced projective morphism 
$Y(\lambda_1) \to Y(\lambda_2)$;
in particular, every $Y(\lambda)$ is projective over $Y(0)$.
\end{enumerate}
\end{theorem}

\index{GIT-fan!}%
The collection  $\Lambda(X,H)$ 
is called the {\em GIT-(quasi-)fan\/}
of the $H$-variety $X$.
For a diagonal $H$-action on $\KK^r$, 
the orbit cones are 
$\cone(\deg(T_i); \; i \in I)$,
where $I$ runs through the subsets of 
$\{1, \ldots, r\}$
and, thus the GIT-fan equals
the Gelfand-Kapranov-Zelevinsky 
decomposition associated to 
$\deg(T_1), \ldots, \deg(T_r) \in K_\QQ$.

\begin{example}
\label{ex:std3onk6}
Consider the action of the standard three torus 
$\TT^3$ on $\KK^6$ defined by $\deg(T_i) = w_i$,
where $w_1, \ldots, w_6$ are defined as
$$ 
w_1 \ = \ (1,0,0),\quad 
w_2 \ = \ (0,1,0),\quad 
w_3 \ = \ (0,0,1),
$$
$$ 
w_4 \ = \ (1,1,0),\quad 
w_5 \ = \ (1,0,1),\quad 
w_6 \ = \ (0,1,1).
$$
Then the GIT-fan subdivides the positive orthant 
in $\QQ^3$; intersecting with a suitable plane 
perpendicular to the line through $(1,1,1)$ 
gives the following picture.
$$
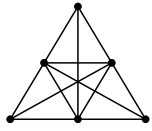
\begin{picture}(0,0)%
\includegraphics{13chambers.pstex}%
\end{picture}%
\setlength{\unitlength}{1243sp}%
\begingroup\makeatletter\ifx\SetFigFontNFSS\undefined%
\gdef\SetFigFontNFSS#1#2#3#4#5{%
  \reset@font\fontsize{#1}{#2pt}%
  \fontfamily{#3}\fontseries{#4}\fontshape{#5}%
  \selectfont}%
\fi\endgroup%
\begin{picture}(2850,2400)(1276,-2836)
\end{picture}%

$$
\end{example}

If every $H$-invariant divisor on $X$ is principal,
then the GIT-fan controls the whole variation of
good quotients with a quasiprojective quotient
space. 
For the precise statement let us call 
a good $H$-set $U \subseteq X$
{\em qp-maximal\/} if 
$U \quot H$ is quasiprojective and 
$U$ is maximal w.r.t. $H$-saturated
inclusion among all good $H$-sets
$W \subseteq X$ with $W \quot H$
quasiprojective.

\begin{theorem}
\label{thm:qpmax}
Assume that $X$ is normal and 
for every $H$-invariant divisor
on~$X$ some positive multiple is principal.
Then, with the GIT-fan $\Lambda(X,H)$ of 
the $H$-action on $X$, we have mutually
inverse order reversing bijections
\begin{eqnarray*}
\Lambda(X,H)  
& \longleftrightarrow &
\{\text{qp-maximal subsets of } X \}
\\
\lambda
& \mapsto & 
X^{ss}(\lambda) = \{x \in X; \; \lambda \subseteq  \omega_x\}
\\
\bigcap_{x \in U} \omega_x =: \lambda(U)
& \mapsfrom &
U.
\end{eqnarray*}
\end{theorem}

Now we look for more general quotient spaces.
We say that a variety $X$ has the
{\em $A_2$-property \/}, if any two points
$x,x' \in X$ admit a common affine open
neighborhood in $X$. 
By~\cite{Wl}, the normal $A_2$-varieties 
are precisely those that admit a closed 
embedding into a toric variety.

\begin{definition}
\label{def:bunchorbcones}
\index{bunch! of orbit cones}%
Let $\Omega_X$ denote the collection of
all orbit cones $\omega_x$, where $x \in X$.
A {\em bunch of orbit cones\/} is 
a nonempty collection $\Phi \subseteq \Omega_X$ 
such that
\begin{enumerate}
\item
given $\omega_1, \omega_2 \in \Phi$, one has
$\omega_1^\circ \cap \omega_2^\circ \ne \emptyset$,
\item
given $\omega \in \Phi$, every orbit cone
$\omega_0 \in \Omega_X$ with 
$\omega^\circ \subseteq  \omega_0^\circ$ 
belongs to $\Phi$.
\end{enumerate}
A {\em maximal bunch of orbit cones\/} is a 
bunch of orbit cones $\Phi \subseteq \Omega_X$ 
which cannot be enlarged by adding further 
orbit cones.
\end{definition}

\begin{definition}
\index{refinement} 
Let $\Phi,\Phi' \subseteq \Omega_X$ be 
bunches of orbit cones.
We say that 
$\Phi$ {\em refines\/} $\Phi'$
(written $\Phi \le \Phi')$, 
if for any $\omega' \in \Phi'$ 
there is an $\omega \in \Phi$
with $\omega \subseteq \omega'$. 
\end{definition}

\begin{example}
\label{ex:std3onk6bunches}
Consider once more the action of the standard 
three torus $\TT^3$ 
on $\KK^6$ discussed in~\ref{ex:std3onk6}.
Here are two maximal bunches, indicated by drawing 
their minimal members:
$$
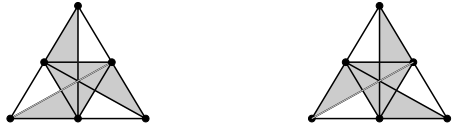
\begin{picture}(0,0)%
\includegraphics{14bunches1.pstex}%
\end{picture}%
\setlength{\unitlength}{1243sp}%
\begingroup\makeatletter\ifx\SetFigFontNFSS\undefined%
\gdef\SetFigFontNFSS#1#2#3#4#5{%
  \reset@font\fontsize{#1}{#2pt}%
  \fontfamily{#3}\fontseries{#4}\fontshape{#5}%
  \selectfont}%
\fi\endgroup%
\begin{picture}(2850,2400)(1276,-2836)
\end{picture}%

\quad\qquad\quad\qquad
\begin{picture}(0,0)%
\includegraphics{14bunches2.pstex}%
\end{picture}%
\setlength{\unitlength}{1243sp}%
\begingroup\makeatletter\ifx\SetFigFontNFSS\undefined%
\gdef\SetFigFontNFSS#1#2#3#4#5{%
  \reset@font\fontsize{#1}{#2pt}%
  \fontfamily{#3}\fontseries{#4}\fontshape{#5}%
  \selectfont}%
\fi\endgroup%
\begin{picture}(2850,2400)(1276,-2836)
\end{picture}%

$$
\end{example}

\begin{definition}
\label{def:UofPsi}
To any collection of orbit cones
$\Phi$ of $X$,
we associate the following 
subset of $X$:
\begin{eqnarray*}
U(\Phi)
& := &
\{ x \in X; \; \omega_0 \preceq \omega_x 
\text{ for some } \omega_0 \in \Phi \}.
\end{eqnarray*}
Conversely, to any $H$-invariant
subset $U \subseteq X$,
we associate the following  collection
of orbit cones
\begin{eqnarray*}
\Phi(U)
& := &
\{\omega_x; \; x \in U \text{ with }
H \mal x \text{ closed in } U
\}.
\end{eqnarray*}
\end{definition}

By an {\em $(H,2)$-maximal subset\/} of $X$
we mean a good $H$-set $U \subseteq X$
with $U \quot H$ an $A_2$-variety
such that $U$ is maximal w.r.t. $H$-saturated
inclusion among all good $H$-sets
$W \subseteq X$ with $W \quot H$
an $A_2$-variety.
We are ready for the result.

\begin{theorem}
\label{thm:maxcoll}
Assume that $X$ is normal and 
for every $H$-invariant divisor
on~$X$ some positive multiple is principal.
Then we have mutually inverse order reversing
bijections
\begin{eqnarray*}
\left\{
 \text{maximal bunches of orbit cones in } \Omega_X
\right\}
& \longleftrightarrow &
\{\text{$(H,2)$-maximal subsets of } X \}
\\
\Phi 
& \mapsto & 
U(\Phi)
\\
\Phi(U)
& \mapsfrom &
U.
\end{eqnarray*}
\end{theorem}

\begin{remark}
\label{rem:chmber2bunch}
Every GIT-chamber $\lambda \in \Lambda(X,H)$
defines a bunch of orbit cones 
$\Phi(\lambda) = \{\omega_x; \; \lambda^\circ \subseteq \omega^\circ\}$.
These bunches turn out to be maximal
and they correspond to the qp-maximal
subsets of $X$; 
in particular, the latter ones
are $(H,2)$-maximal.
The bunches of~\ref{ex:std3onk6bunches} give 
rise to non-projective complete quotients.
\end{remark}

\subsection{Cox rings and combinatorics}
Here we present the combinatorial approach 
to varieties with finitely generated Cox 
ring developed in~\cite{BeHa4,Ha2},
see also~\cite[Chap.~III]{ArDeHaLa}.
The approach generalizes the combinatorial description 
of toric varieties and has many common features 
with methods of~\cite{Do,IaFl} for investigating 
subvarieties of weighted complete spaces.
The whole thing is based on the following simple
observation.

\begin{remark}
Let $X$ be a normal variety with 
$\Gamma(X,\mathcal{O}^*) = \KK^*$
and finitely generated divisor class 
group. 
If the Cox ring $\mathcal{R}(X)$ is 
finitely generated, then we obtain
the following picture
$$ 
\xymatrix{
{\Spec_X \mathcal{R}}
\ar@{}[r]|{\quad=}
& 
{\rq{X}}
\ar@{}[r]|\subseteq
\ar[d]_{\quot H_X}
&
{\b{X}}
\ar@{}[r]|{=\quad}
&
{\Spec \, \mathcal{R}(X)}
\\
& 
X
&
&
}
$$
where $\rq{X} \subseteq \b{X}$ is an open 
$H_X$-invariant subset of the $H_X$-factorial
affine variety $\b{X}$ and the characteristic 
space $\rq{X} \to X$ is a good quotient for
the $H_X$-action. We call the affine 
$H_X$-variety $\b{X}$ the {\em total coordinate space\/}
of $X$. 
\end{remark}

Thus, we see that all varieties sharing the same 
divisor class group $K$ and finitely generated Cox 
ring $R$ occur as good quotients of suitable open subsets of 
$\Spec \, R$ by the action of $\Spec \, \KK[K]$.
The latter ones we just described in combinatorial
terms via Geometric Invariant Theory.
We now turn this picture into a combinatorial 
language allowing explicit computations.

\begin{definition}
\label{def:brdata}
Let $K$ be a finitely generated abelian group
and $R$ a factorially $K$-graded affine algebra 
with $R^* = \KK^*$.
Moreover, let $\mathfrak{F} = (f_1, \ldots, f_r)$ 
be a system of pairwise nonassociated 
$K$-prime generators for $R$.
\begin{enumerate}
\item 
\index{projected cone!}%
\index{cone! projected}%
The {\em projected cone\/} associated
to $\mathfrak{F}$ is 
$(E \topto{Q} K, \gamma)$,
where $ E := \ZZ^{r}$, the homomorphism
$Q \colon E \to K$ sends the $i$-th canonical
basis vector $e_{i} \in E$ to  
$w_i := \deg(f_{i}) \in K$
and $\gamma \subseteq  E_{\QQ}$ is the 
convex cone generated by $e_{1}, \dots, e_{r}$.
\item
\index{almost free grading!}%
\index{grading! almost free}
We say that the $K$-grading of $R$ is 
{\em almost free\/} if  for every facet 
$\gamma_0 \preceq \gamma$ the image 
$Q(\gamma_0 \cap E)$ generates the 
abelian group~$K$.
\item
We say that $\gamma_0 \preceq \gamma$ is an 
{\em $\mathfrak{F}$-face}, if the product 
over all $f_{i}$ with $e_{i} \in \gamma_{0}$
does not lie in 
$\sqrt{\bangle{f_{j}; \; e_{j} \not\in \gamma_{0}}} \subseteq A$.
\item
\index{$\mathfrak{F}$-bunch!}%
Let 
$\Omega_{\mathfrak{F}} 
= 
\{Q(\gamma_0); \; 
\gamma_0 \preceq \gamma \ \mathfrak{F} \text{-face}\}$
denote the collection of projected $\mathfrak{F}$-faces.
An {\em $\mathfrak{F}$-bunch\/} is 
a nonempty subset $\Phi \subseteq \Omega_{\mathfrak{F}}$ 
such that
\begin{enumerate}
\item
for any two $\tau_1,\tau_2 \in \Phi$, we have 
$\tau_1^{\circ} \cap \tau_2^{\circ} \ne \emptyset$,
\item
if $\tau_1^{\circ} \subseteq \tau^{\circ}$ holds
for $\tau_1 \in \Phi$ and $\tau \in \Omega_{\mathfrak{F}}$,
then $\tau \in \Phi$ holds.
\end{enumerate}
\item
\index{true $\mathfrak{F}$-bunch!}%
\index{$\mathfrak{F}$-bunch! true}%
We say that an $\mathfrak{F}$-bunch
$\Phi$ is {\em true\/}
if for every facet $\gamma_{0} \preceq \gamma$
the image $Q(\gamma_{0})$ belongs to $\Phi$. 
\end{enumerate}
\end{definition}

\begin{definition}
\index{bunched ring!}
\index{ring! bunched}
A {\em bunched ring\/} is a triple 
$(R,\mathfrak{F},\Phi)$, where 
$R$ is an almost freely factorially $K$-graded 
affine $\KK$-algebra such that $R^* = \KK^*$ 
holds,
$\mathfrak{F}$ is a system of pairwise 
non-associated $K$-prime generators 
for $R$ and $\Phi$ is a true 
$\mathfrak{F}$-bunch.
\end{definition}

\begin{construction}
\label{bringvarconstr}
\index{relevant faces!}%
\index{covering collection!}%
Let $(R,\mathfrak{F},\Phi)$ be a bunched ring. 
Then $\Phi$ is a bunch of orbit cones
for the action of $H := \Spec \, \KK[K]$ 
on $\b{X} := \Spec \, R$.
Thus, we have the associated open set and its 
quotient
$$ 
\rq{X} 
\ := \ 
\rq{X}(R,\mathfrak{F},\Phi) 
\ = \ 
\b{X}(\Phi)
\ \subseteq \ 
\b{X},
$$
$$
X
\ := \ 
X(R,\mathfrak{F},\Phi)
\ := \
\rq{X}(R,\mathfrak{F},\Phi) \quot H.
$$
We denote the quotient map by $p \colon \rq{X} \to X$.
Conditions~\ref{def:brdata}~(ii) and~(v) ensure 
that the $H$-action on $\rq{X}$ is strongly stable.
Moreover, every member $f_i$ of $\mathfrak{F}$ defines 
a prime divisor $D_X^i := p(V(\rq{X},f_i))$ on $X$.  
\end{construction}

\begin{theorem}
\label{thm:bring2var}
Let $\rq{X} := \rq{X}(R,\mathfrak{F},\Phi)$
and 
$X := X(R,\mathfrak{F},\Phi)$
arise from a bunched ring
$(R,\mathfrak{F},\Phi)$.
Then $X$ is a normal $A_2$-variety with 
$$
\dim(X) \ = \ \dim(R) - \dim(K_\QQ),
\qquad\qquad
\Gamma(X,\mathcal{O}^*) = \KK^*,
$$
there is an isomorphism 
$\Cl(X) \to K$ sending $[D_X^i]$ to $\deg(f_i)$,
the map $p \colon \rq{X} \to X$ 
is a characteristic space and
the Cox ring $\mathcal{R}(X)$ is isomorphic
to $R$.
\end{theorem}

\begin{theorem}
Every complete normal $A_2$-variety with finitely 
generated Cox ring arises from a bunched ring;
in particular, every projective normal variety 
with finitely generated Cox ring does so.
\end{theorem}

Let us illustrate Construction~\ref{bringvarconstr}
with two examples.
The first one shows how toric varieties
fit into the picture of bunched rings.

\begin{example}[Bunched polynomial rings]
\label{ex:bunchedpolring}
Consider a bunched ring $(R,\mathfrak{F},\Phi)$
with $R = \KK[T_1, \ldots, T_r]$ and 
$\mathfrak{F} := (T_1, \ldots, T_r)$.
Then $X(R,\mathfrak{F},\Phi)$ is a toric 
variety.
Its defining fan $\Sigma$ is obtained from 
$\Phi$ via linear Gale duality:
$$ 
\xymatrix@R=10pt{
0 \ar[rr]
&&
L_\QQ
\ar[rr]^{Q^*}
&&
F_\QQ
\ar[rr]^{P}
&&
N_\QQ
\ar[rr]
&&
0
\\
&&
&&
{\Sigma^{\scriptscriptstyle \uparrow}}
\ar[rr]^{\delta_0 \mapsto P(\delta_0)}
&&
{\Sigma}
&&
\\
\\
&&
{\Phi}
\ar@{<-}[rr]_{Q(\gamma_0) \mapsfrom \gamma_0}
&&
{\Phi^{\scriptscriptstyle \uparrow}}
\ar[uu]_{\gamma_0 \mapsto \gamma_0^\perp \cap \delta}
&&
&&
\\
0 \ar@{<-}[rr]
&&
K_\QQ
\ar@{<-}[rr]_{Q}
&&
E_\QQ
\ar@{<-}[rr]_{P^*}
&&
M_\QQ
\ar@{<-}[rr]
&&
0
}
$$
Here $\Phi^{\scriptscriptstyle \uparrow}$ 
consists of those faces 
of the orthant $\gamma \subseteq E_\QQ$
that map onto a member of $\Phi$
and $\Sigma^{\scriptscriptstyle \uparrow}$ 
of the corresponding faces of the dual orthant 
$\delta \subseteq F_{\QQ}$.
Note that $P \colon F \to N$ is the same map 
as in~\ref{constr:toriccox} and the 
$D_X^i$ are exactly the toric prime divisors.
\end{example}

\begin{example}[A singular del Pezzo surface]
\label{ex:delpezzo}
Consider $K := \ZZ^2$ and the $K$-grading of 
$\KK[T_1, \ldots, T_5]$ given by $\deg(T_i) := w_i$, 
where $w_i$ is the $i$-th column of 
\begin{eqnarray*}
Q
& := & 
\left[
\begin{array}{rrrrr}
1 & -1 & 0 & -1 & 1
\\
1 & 1 & 1 & 0 & 2
\end{array}
\right]
\end{eqnarray*}
Then this $K$-grading descends to a $K$-grading of 
the following residue algebra which is known to 
be factorial:
\begin{eqnarray*}
R 
& := & 
\KK[T_1, \ldots, T_5] \, / \, \bangle{T_1T_2 + T_3^2 + T_4T_5}.
\end{eqnarray*}
The classes $f_i \in R$ of $T_i \in \KK[T_1, \ldots, T_5]$,
where $1 \le i \le 5$, form a system $\mathfrak{F}$ 
of pairwise nonassociated $K$-prime generators
of $R$. We have 
$$ 
E \ = \ \ZZ^5, \qquad \gamma \ = \ \cone(e_1, \ldots, e_5)
$$
and the $K$-grading is almost free. 
Computing the $\mathfrak{F}$-faces,
we see that there is one maximal true
$\mathfrak{F}$-bunch $\Phi$; 
it has $\tau :=  \cone(w_2,w_5)$ as 
its unique minimal cone.
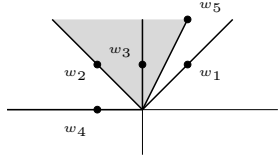
\begin{center}
\begin{picture}(0,0)%
\includegraphics{delpezzobunch.pstex}%
\end{picture}%
\setlength{\unitlength}{1243sp}%
\begingroup\makeatletter\ifx\SetFigFontNFSS\undefined%
\gdef\SetFigFontNFSS#1#2#3#4#5{%
  \reset@font\fontsize{#1}{#2pt}%
  \fontfamily{#3}\fontseries{#4}\fontshape{#5}%
  \selectfont}%
\fi\endgroup%
\begin{picture}(5445,3144)(2218,-2323)
\put(6076,-736){\makebox(0,0)[lb]{\smash{{\SetFigFontNFSS{5}{6.0}{\familydefault}{\mddefault}{\updefault}{\color[rgb]{0,0,0}$w_1$}%
}}}}
\put(6076,614){\makebox(0,0)[lb]{\smash{{\SetFigFontNFSS{5}{6.0}{\familydefault}{\mddefault}{\updefault}{\color[rgb]{0,0,0}$w_5$}%
}}}}
\put(3376,-1861){\makebox(0,0)[lb]{\smash{{\SetFigFontNFSS{5}{6.0}{\familydefault}{\mddefault}{\updefault}{\color[rgb]{0,0,0}$w_4$}%
}}}}
\put(4276,-286){\makebox(0,0)[lb]{\smash{{\SetFigFontNFSS{5}{6.0}{\familydefault}{\mddefault}{\updefault}{\color[rgb]{0,0,0}$w_3$}%
}}}}
\put(3376,-736){\makebox(0,0)[lb]{\smash{{\SetFigFontNFSS{5}{6.0}{\familydefault}{\mddefault}{\updefault}{\color[rgb]{0,0,0}$w_2$}%
}}}}
\end{picture}%

\end{center}
Note that $\tau$ is a GIT-cone,
$\Phi=\Phi(\tau)$ holds with $\Phi(\tau)$ 
as in~\ref{rem:chmber2bunch}
and $\rq{X}(R,\mathfrak{F},\Phi)$ 
equals $\b{X}^{ss}(\tau)$ in 
$\b{X} = V(\KK^5; T_1T_2 + T_3^2 + T_4T_5)$.
For $X = X(R,\mathfrak{F},\Phi)$
we have
$$ 
\dim(X) \ = \ 2,
\qquad
\Cl(X) \ = \ \ZZ^2,
\qquad
\mathcal{R}(X) \ = \ R.
$$
\end{example}

\begin{definition}
Let $(R,\mathfrak{F},\Phi)$ be a bunched ring 
and $(E \topto{Q} K,\gamma)$ its projected cone.
The {\em collection of relevant faces\/} 
and the {\em covering collection\/} are
\begin{eqnarray*}
\rel(\Phi)
& := & 
\{
\gamma_0 \preceq \gamma; \; 
\gamma_0 \text{ an $\mathfrak{F}$-face with }
Q(\gamma_0) \in \Phi
\},
\\
\cov(\Phi)
& := & 
\{
\gamma_0 \in \rel(\Phi); \; 
\gamma_0 \text{ minimal} \}.
\end{eqnarray*}
\end{definition}

\begin{construction}[Canonical toric embedding]
Any bunched ring $(R,\mathfrak{F},\Phi)$ defines 
a bunched polynomial ring $(R',\mathfrak{F}',\Phi')$
by ``forgetting the relations'':
If the system of generators of $R$ is 
$\mathfrak{F} = (f_1, \ldots, f_r)$, set
$$
R' := \KK[T_1,\ldots, T_r],
\qquad
\deg(T_i) := \deg(f_i) \in K,
\qquad
\mathfrak{F}' := (T_1, \ldots, T_r)
$$ 
and let $\Phi'$ be the $\mathfrak{F}'$-bunch 
generated by $\Phi$, i.e. it consists of all
projected faces $Q(\gamma_0)$ with 
$\tau^\circ \subseteq Q(\gamma_0)^\circ$
for some $\tau \in \Phi$.
Then we obtain a commutative diagram,
where the induced map of quotients 
$\imath \colon X \to Z$ is a closed 
embedding of the varieties $X$ and $Z$ 
associated to the bunched rings 
$(R,\mathfrak{F},\Phi)$ and
$(R',\mathfrak{F}',\Phi')$ respectively:
$$ 
\xymatrix{
{\b{X}}
\ar@{}[r]|\supseteq
&
{\rq{X}}
\ar[r]
\ar[d]_{\quot H}
& 
{\rq{Z}}
\ar[d]^{\quot H}
\ar@{}[r]|\subseteq
&
{\b{Z}}
\\
&
X 
\ar[r]_{\imath}
&
Z
&
}
$$
By construction, $Z$ is toric and 
we have an isomorphism 
$\imath^* \colon \Cl(Z) \to \Cl(X)$.
Then, in the setting of~\ref{ex:bunchedpolring},
the toric orbits of $Z$ 
intersecting $X$ are precisely the orbits
$B(\sigma) \subseteq Z$ corresponding to cones 
$\sigma = P(\gamma_0^*)$ with 
$\gamma_0 \in \rel(\Phi)$.
In particular, we obtain a decomposition
into locally closed strata
$$ 
X 
\ = \ 
\bigcup_{\gamma_0 \in \rel(\Phi)} X(\gamma_0),
\qquad
X(\gamma_0) := X \cap B(\sigma).
$$
\end{construction}

\begin{remark}
In general, the canonical toric ambient 
variety $Z$ is not complete,
even if $X$ is.
If $X$ is projective, then $\Phi = \Phi(\lambda)$
holds with a GIT-cone $\lambda \in \Lambda(\b{X},H)$.
The toric GIT-fan $\Lambda(\b{Z},H)$ refines 
$\Lambda(\b{X},H)$ and every 
$\eta \in \Lambda(\b{Z},H)$ with 
$\eta^\circ \subseteq \lambda^\circ$
defines a projective completion of $Z$.
For example, in the setting of~\ref{ex:delpezzo},
the two GIT-fans are
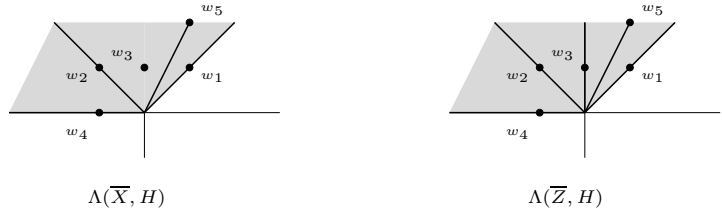
\begin{center}
\begin{picture}(0,0)%
\includegraphics{delptoricgitfan.pstex}%
\end{picture}%
\setlength{\unitlength}{1243sp}%
\begingroup\makeatletter\ifx\SetFigFontNFSS\undefined%
\gdef\SetFigFontNFSS#1#2#3#4#5{%
  \reset@font\fontsize{#1}{#2pt}%
  \fontfamily{#3}\fontseries{#4}\fontshape{#5}%
  \selectfont}%
\fi\endgroup%
\begin{picture}(5445,4159)(2218,-3338)
\put(6076,-736){\makebox(0,0)[lb]{\smash{{\SetFigFontNFSS{5}{6.0}{\familydefault}{\mddefault}{\updefault}{\color[rgb]{0,0,0}$w_1$}%
}}}}
\put(6076,614){\makebox(0,0)[lb]{\smash{{\SetFigFontNFSS{5}{6.0}{\familydefault}{\mddefault}{\updefault}{\color[rgb]{0,0,0}$w_5$}%
}}}}
\put(3376,-1861){\makebox(0,0)[lb]{\smash{{\SetFigFontNFSS{5}{6.0}{\familydefault}{\mddefault}{\updefault}{\color[rgb]{0,0,0}$w_4$}%
}}}}
\put(4276,-286){\makebox(0,0)[lb]{\smash{{\SetFigFontNFSS{5}{6.0}{\familydefault}{\mddefault}{\updefault}{\color[rgb]{0,0,0}$w_3$}%
}}}}
\put(3376,-736){\makebox(0,0)[lb]{\smash{{\SetFigFontNFSS{5}{6.0}{\familydefault}{\mddefault}{\updefault}{\color[rgb]{0,0,0}$w_2$}%
}}}}
\put(3826,-3211){\makebox(0,0)[lb]{\smash{{\SetFigFontNFSS{7}{8.4}{\familydefault}{\mddefault}{\updefault}{\color[rgb]{0,0,0}$\Lambda(\b{X},H)$}%
}}}}
\end{picture}%

\qquad \qquad \qquad
\begin{picture}(0,0)%
\includegraphics{delpezzogitfan.pstex}%
\end{picture}%
\setlength{\unitlength}{1243sp}%
\begingroup\makeatletter\ifx\SetFigFontNFSS\undefined%
\gdef\SetFigFontNFSS#1#2#3#4#5{%
  \reset@font\fontsize{#1}{#2pt}%
  \fontfamily{#3}\fontseries{#4}\fontshape{#5}%
  \selectfont}%
\fi\endgroup%
\begin{picture}(5445,4159)(2218,-3338)
\put(6076,-736){\makebox(0,0)[lb]{\smash{{\SetFigFontNFSS{5}{6.0}{\familydefault}{\mddefault}{\updefault}{\color[rgb]{0,0,0}$w_1$}%
}}}}
\put(6076,614){\makebox(0,0)[lb]{\smash{{\SetFigFontNFSS{5}{6.0}{\familydefault}{\mddefault}{\updefault}{\color[rgb]{0,0,0}$w_5$}%
}}}}
\put(3376,-1861){\makebox(0,0)[lb]{\smash{{\SetFigFontNFSS{5}{6.0}{\familydefault}{\mddefault}{\updefault}{\color[rgb]{0,0,0}$w_4$}%
}}}}
\put(4276,-286){\makebox(0,0)[lb]{\smash{{\SetFigFontNFSS{5}{6.0}{\familydefault}{\mddefault}{\updefault}{\color[rgb]{0,0,0}$w_3$}%
}}}}
\put(3376,-736){\makebox(0,0)[lb]{\smash{{\SetFigFontNFSS{5}{6.0}{\familydefault}{\mddefault}{\updefault}{\color[rgb]{0,0,0}$w_2$}%
}}}}
\put(3826,-3211){\makebox(0,0)[lb]{\smash{{\SetFigFontNFSS{7}{8.4}{\familydefault}{\mddefault}{\updefault}{\color[rgb]{0,0,0}$\Lambda(\b{Z},H)$}%
}}}}
\end{picture}%

\end{center}
The GIT-cones $\cone(w_2,w_3)$ and $\cone(w_3,w_5)$
in $\Lambda(\b{Z},H)$ provide completions of $Z$
by $\QQ$-factorial projective toric
varieties $Z_1$ and $Z_2$ and 
$\cone(w_3)$ gives a completion
by a projective toric variety $Z_3$
with a non-$\QQ$-factorial singularity.
\end{remark}

\goodbreak

We now indicate how to read off
basic geometric properties from
defining data.
In the sequel, $X$ is the variety 
arising from a bunched ring 
$(R,\mathfrak{F},\Phi)$.

\begin{theorem}
\label{thm:localcl}
Consider a relevant face $\gamma_0 \in \rel(\Phi)$ and 
a point $x \in X(\gamma_0)$.
Then we have a commutative diagram
$$ 
\xymatrix{
{\Cl(X)}
\ar[r]
\ar@{<->}[d]_{\cong}
&
{\Cl(X,x)}
\ar@{<->}[d]^{\cong}
\\
K
\ar[r]
&
K / Q(\lin(\gamma_0) \cap E))
}
$$
In particular, the local divisor class groups 
are constant along the pieces $X(\gamma_0)$,
where $\gamma_0 \in \rel(\Phi)$.
Moreover, the Picard group of $X$ is given by
\begin{eqnarray*}
\Pic(X) 
& \cong & 
\bigcap_{\gamma_{0} \in \cov(\Phi)} 
              Q(\lin(\gamma_{0}) \cap E).
\end{eqnarray*}
\end{theorem}

\begin{theorem}
\label{singularities}
Consider a relevant face $\gamma_{0} \in \rel(\Phi)$
and point $x \in X(\gamma_{0})$ in the corresponding 
stratum.
\begin{enumerate}
\item 
The point $x$ is factorial if and only if 
$Q$ maps $\lin(\gamma_{0}) \cap E$ onto $K$.
\item 
The point $x$ is $\QQ$-factorial if and only if 
 $Q(\gamma_{0})$ is of full dimension.
\end{enumerate}
In particular, $X$ is $\QQ$-factorial, if and 
only if $\Phi$ consists of full-dimensional cones.
If $\rq{X}$ is smooth, then every factorial
point of $X$ is smooth.
\end{theorem}

\begin{theorem}
\label{thm:ample}
In the divisor class group $K = \Cl(X)$, we have the
following descriptions of the cones 
of effective, movable, semiample and ample divisors:
$$
\begin{array}{ccc}
\displaystyle \Eff(X) \ = \ Q(\gamma),
& \qquad &
\displaystyle \Mov(X) 
\; = 
\bigcap_{\gamma_0 \text{\ facet of \ } \gamma} Q(\gamma_0),
\\
 & & 
\\
\displaystyle \SAmple(X) \ = \ \bigcap_{\tau \in \Phi} \tau,
& \qquad &
\displaystyle \Ample(X) \ = \ \bigcap_{\tau \in \Phi} \tau^{\circ}.
\end{array}
$$
\end{theorem}

\begin{example}[The singular del Pezzo surface, continued]
The variety $X$ of~\ref{ex:delpezzo} is a $\QQ$-factorial
surface.
It has a single singularity, namely the point in the 
piece $x_0 \in X(\gamma_0)$ for $\gamma_0 = \cone(e_2,e_5)$.
The local class group $\Cl(X,x_0)$ is cyclic of order three
and the Picard group of $X$ is of index $3$ in $\Cl(X)$.
Moreover, the ample cone of $X$ is generated by $w_2$ and 
$w_5$. In particular, $X$ is projective.
\end{example}

\begin{remark}
Applying~\ref{thm:localcl}, \ref{singularities} 
and~\ref{thm:ample} to the canonical toric ambient 
variety $Z$ shows that $X$ inherits local class groups 
and singularities from $Z$ and the Picard group as 
well as the various cones of divisors of $X$ and 
$Z$ coincide.
However, factorial singularities of $X$ 
are smooth points of $Z$, see~\ref{ex:E6cubicresol}
for an example.
Moreover, $\Pic(X) = \Pic(Z)$ and 
$\Ample(X) = \Ample(Z)$ can get lost when replacing 
$Z$ with a completion. 
\end{remark}

The following two statements concern the 
case that $R$ is a complete intersection
in the sense that 
with $d := r - \dim(X) - \dim(K)$, 
there are $K$-homogeneous generators 
$g_{1}, \ldots, g_{d}$ for the ideal of 
relations between $f_{1}, \ldots, f_{r}$. 
Set $w_i := \deg(f_i)$ and $u_j := \deg(g_j)$.

\begin{theorem}
\label{thm:canondiv}
Suppose that $R$ is a complete intersection as 
above. Then the canonical class of $X$ is 
given in $K = \Cl(X)$ by
\begin{eqnarray*}
\mathcal{K}_X
& = &
\sum u_j - \sum w_i.
\end{eqnarray*}
\end{theorem}

\goodbreak

\begin{remark}[Computing intersection numbers]
\label{rem:compint}
Suppose that $R$ is a complete intersection 
and that $\Phi = \Phi(\lambda)$ 
holds with a full-dimensional  
$\lambda \in \Lambda(\b{X},H)$.
Fix a full-dimensional $\eta \in \Lambda(\b{Z},H)$ 
with $\eta^\circ \subseteq \lambda^\circ$.
For $w_{i_1}, \ldots, w_{i_{n+d}}$
let $w_{j_1}, \ldots, w_{j_{r-n-d}}$ denote
the complementary weights and set
\begin{eqnarray*}
\tau(w_{i_1}, \ldots, w_{i_{n+d}})
& := & 
\cone(w_{j_1}, \ldots, w_{j_{r-n-d}}),
\\
\mu(w_{i_1}, \ldots, w_{i_{n+d}})
& := & 
[K : \bangle{w_{j_1}, \ldots, w_{j_{r-n-d}}}].
\end{eqnarray*}
Then the intersection product $K_\QQ^{n+d} \to \QQ$
of the ($\QQ$-factorial) toric variety $Z_1$ 
associated to $\Phi(\eta)$ is determined by 
the values
\begin{eqnarray*}
w_{i_1} \cdots w_{i_{n+d}}
& = & 
\begin{cases}
\mu(w_{i_1}, \ldots, w_{i_{n+d}})^{-1},
& 
\eta \subseteq \tau(w_{i_1}, \ldots, w_{i_{n+d}}),
\\
0,
& 
\eta \not\subseteq \tau(w_{i_1}, \ldots, w_{i_{n+d}}).
\end{cases}
\end{eqnarray*}
As a complete intersection, $X \subseteq Z_1$
inherits intersection theory.  
For a tuple $D_X^{i_1}, \ldots, D_X^{i_n}$
on $X$, its intersection 
number can be computed by 
\begin{eqnarray*}
D_X^{i_1} \cdots D_X^{i_n} 
& = & 
w_{i_1} \cdots w_{i_n} \cdot u_1 \cdots u_d.
\end{eqnarray*}
\end{remark}

\begin{example}[The singular del Pezzo surface, continued]
\label{ex:delpezzo2}
Consider once more the surface $X$ of~\ref{ex:delpezzo}.
The degree of the defining relation is 
$\deg(T_1T_2 + T_3^2 + T_4T_5) = 2w_3$
and thus the canonical class of $X$ 
is given as 
$$ 
\mathcal{K}_X
\ = \
2w_3
-
(w_1 + w_2 + w_3 + w_4 + w_5)
\ = \ 
- 3w_3.
$$
In particular, we see that the anticanonical class is ample
and thus $X$ is a (singular) del Pezzo surface.
The self intersection number of the canonical class is
$$ 
\mathcal{K}_X^2
\ = \ 
(3w_3)^2
\ = \ 
\frac{9(w_1+w_2)(w_4+w_5)}{4}.
$$
The $w_i \cdot w_j$ equal the toric 
intersection numbers 
$2 w_i \cdot w_j \cdot w_3$. 
In order to compute them,
let $w^1_{ij}, w^2_{ij}$ denote the weights
in $\{w_1, \ldots, w_5\} \setminus \{w_i,w_j,w_3\}$.
Then we have
\begin{eqnarray*}
w_i \cdot w_j \cdot w_3
& = & 
\begin{cases}
\mu(w_i,w_j,w_3)^{-1}, &  \tau \subseteq \cone (w^1_{ij}, w^2_{ij}),
\\
0, & \tau \not\subseteq \cone (w^1_{ij}, w^2_{ij}), 
\end{cases}
\end{eqnarray*}
where the multiplicity $\mu(w_i,w_j,w_3)$ is 
the absolute value of $\det(w^1_{ij}, w^2_{ij})$.
Thus, we can proceed in the computation:
$$ 
\mathcal{K}_X^2
\ =  \
(3w_3)^2
\ = \ 
\frac{9 \cdot 2}{4} 
\left( 
\vert \det(w_2,w_5) \vert^{-1} 
+  
\vert \det(w_1,w_4) \vert^{-1}  
\right)
\ = \ 
\frac{9}{2} \cdot \frac{4}{3}
\ = \ 
6.
$$
\end{example}


\subsection{Mori dream spaces}
We take a closer look to the $\QQ$-factorial 
projective varieties with a finitely generated 
Cox ring.
Hu and Keel~\cite{HuKe} called them {\em Mori dream spaces\/} 
and characterized them in terms of cones of divisors.
In the context of normal complete varieties 
their statement reads as follows; 
by a {\em small birational map\/} 
we mean a rational map defining an isomorphism of open
subsets with complement of codimension two.

\begin{theorem}
\label{fingenchar}
Let $X$ be a normal complete variety 
with finitely generated divisor class group. 
Then the following statements are equivalent.
\begin{enumerate}
\item
The Cox ring $\mathcal{R}(X)$ is finitely generated.
\item
There are small birational maps $\pi_i \colon X \to X_i$, 
where $i = 1, \ldots, r$,
such that each semiample cone 
$\SAmple(X_i) \subseteq \Cl_{\QQ}(X)$ is polyhedral and 
one has
\begin{eqnarray*}
\Mov(X)
& = & 
\pi_1^*(\SAmple(X_1)) \ \cup \ \ldots \ \cup \ \pi_r^*(\SAmple(X_r)).
\end{eqnarray*}
\end{enumerate}
Moreover, if one of these two statements holds, 
then there is a small birational map $X \to X'$
with a $\QQ$-factorial projective variety $X'$.
\end{theorem}

There are many examples of Mori dream spaces. 
Besides the toric and more generally spherical 
varieties, all unirational varieties with a 
reductive group action of complexity one are 
Mori dream spaces.
Other important examples are the log terminal 
Fano varieties~\cite{BCHM}.
Moreover, K3- and Enriques surfaces are 
Mori dream spaces if and only if their effective 
cone is polyhedral~\cite{ArHaLa,To};
the same holds in higher dimensions for  
Calabi-Yau varieties~\cite{McK}.
Specializing to surfaces, the above theorem 
gives the following.

\begin{corollary}
\label{normsurffingen}
Let $X$ be a normal complete surface with 
finitely generated divisor class group 
$\Cl(X)$. 
Then the following statements are
equivalent.
\begin{enumerate}
\item
The Cox ring $\mathcal{R}(X)$ is finitely generated.
\item 
One has  $\Mov(X) = \SAmple(X)$ and this cone is polyhedral.
\end{enumerate}
Moreover, if one of these two statements holds, then 
the surface $X$ is $\QQ$-factorial and projective.
\end{corollary}

The Mori dream spaces sharing a given Cox ring fit 
into a nice picture in terms of the GIT-fan; 
by the {\em moving cone\/} of the $K$-graded 
algebra $R$ we mean here the intersection $\Mov(R)$ 
over all 
$\cone(w_1, \ldots, \rq{w_i}, \ldots, w_r)$,
where the $w_i$ are the degrees of any system of 
pairwise nonassociated homogeneous $K$-prime
generators for $R$.

\begin{remark}
Let $R = \oplus_K R_w$ be an almost freely 
factorially graded affine algebra with 
$R_0 = \KK$ and consider the GIT-fan 
$\Lambda(\b{X},H)$ of the action of 
$H = \Spec \, \KK[K]$ on $\b{X} = \Spec \, R$.
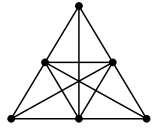
\begin{center}
\begin{picture}(0,0)%
\includegraphics{13chambers.pstex}%
\end{picture}%
\setlength{\unitlength}{1243sp}%
\begingroup\makeatletter\ifx\SetFigFontNFSS\undefined%
\gdef\SetFigFontNFSS#1#2#3#4#5{%
  \reset@font\fontsize{#1}{#2pt}%
  \fontfamily{#3}\fontseries{#4}\fontshape{#5}%
  \selectfont}%
\fi\endgroup%
\begin{picture}(2850,2400)(1276,-2836)
\end{picture}%

\end{center}
Then every GIT-cone $\lambda \in \Lambda(\b{X},H)$ 
defines a projective variety 
$X(\lambda) = \b{X}^{ss}(\lambda) \quot H$.
If $\lambda^\circ \subseteq \Mov(R)^\circ$ holds,
then $X(\lambda)$ is the variety associated to the 
bunched ring $(R,\mathfrak{F},\Phi(\lambda))$ with 
$\Phi(\lambda)$ defined as in~\ref{rem:chmber2bunch}. 
In particular, in this case we have
$$ 
\Cl(X(\lambda)) \ = \ K,
\qquad 
\mathcal{R}(X(\lambda)) \ =  \ R,
$$
$$
\Mov(X(\lambda)) \ = \ \Mov(R),
\qquad
\SAmple(X(\lambda)) \ = \ \lambda.
$$
All projective varieties with Cox ring $R$ are 
isomorphic to some $X(\lambda)$
with $\lambda^\circ \subseteq \Mov(R)^\circ$
and the Mori dream spaces among them 
are precisely those arising from a full dimensional 
$\lambda$.
\end{remark}

Let $X$ be the variety arising from a 
bunched ring $(R,\mathfrak{F},\Phi)$.
Every Weil divisor $D$ on $X$ defines 
a positively graded sheaf
$$ 
\mathcal{S}^+
\ := \ 
\bigoplus_{n \in \ZZ_{\ge 0}} \mathcal{S}^+_n,
\qquad\qquad
\mathcal{S}^+_n
\ := \ 
\mathcal{O}_{X}(nD).
$$
The algebra of global sections of this sheaf inherits 
finite generation from the Cox ring.
In particular, we obtain a rational map
$$ 
\varphi(D) \colon X  \ \to \ X(D),
\qquad
X(D) \ := \ \Proj(\Gamma(X,\mathcal{S}^+)).
$$
Note that $X(D)$ is explicitly given as 
the closure of the image of 
the rational map $X \to \PP_m$ 
determined by the linear system of a 
sufficiently big multiple $nD$.

\begin{remark}
Consider the GIT-fan $\Lambda(\b{X},H)$ 
of the action of $H = \Spec \, \KK[K]$ 
on $\b{X} = \Spec \, R$. Let 
$\lambda  \in \Lambda(\b{X},H)$ be the cone 
with $[D] \in \lambda^{\circ}$ and
$W \subseteq \b{X}$ the open subset 
obtained by removing the zero sets of the 
generators $f_1, \ldots, f_r \in R$.
Then we obtain a commutative diagram
$$ 
\xymatrix{
{\rq{X}}
\ar@{}[r]|\supseteq
\ar[d]_{\quot H}
&
W
\ar@{}[r]|\subseteq
\ar[d]
&
{\b{X}^{ss}(\lambda)}
\ar[d]^{\quot H}
\\
X
\ar@{}[r]|\supseteq
\ar@/_1pc/@{-->}[drr]_{\varphi(D)}
&
W/H
\ar[r]
&
X(\lambda)
\ar@{=}[d]
&
\\
& 
& 
X(D)
}
$$
\end{remark}

\begin{proposition}
Let $D \in \WDiv(X)$ be any Weil divisor,
and denote by $[D] \in \Cl(X)$ its class. 
Then the associated rational map 
$\varphi(D) \colon X \to X(D)$ is 
\begin{enumerate}
\item
birational if and only if 
$[D] \in \Eff(X)^{\circ}$ holds,
\item
small birational if and only if 
$[D] \in \Mov(X)^{\circ}$ holds,
\item
a morphism if and only if 
$[D] \in \SAmple(X)$ holds,
\item
an isomorphism if and only if 
$[D] \in \Ample(X)$ holds.
\end{enumerate}
\end{proposition}

\begin{remark}
The observations made so far imply in particular
that two Mori dream surfaces are isomorphic
if and only if their Cox rings are isomorphic
as graded rings.
Moreover they prove the implication 
``(i)$\Rightarrow$(ii)'' of Theorem~\ref{fingenchar}.
For the other direction, one reduces
finite generation of the Cox ring $\mathcal{R}(X)$ 
to finite generation of the semiample subalgebras 
$\oplus_{K \cap \lambda} \Gamma(X,\mathcal{R}_{[D]})$, 
where $\lambda^\circ \subseteq \Mov(X)^\circ$,
which is given by classical results.
\end{remark}

Recall that two Weil divisors $D,D' \in \WDiv(X)$ 
are said to be {\em Mori equivalent\/}, 
if there is a commutative diagram 
$$ 
\xymatrix{
& 
X 
\ar@{-->}[dl]_{\varphi(D)}
\ar@{-->}[dr]^{\varphi(D')}
& 
\\
X(D)
\ar@{<->}[rr]_{\cong}
&
&
X(D')
}
$$

\begin{proposition}
For any two Weil divisors $D,D'$ on $X$,
the following statements are 
equivalent.
\begin{enumerate}
\item
The divisors $D$ and $D'$ are Mori equivalent.
\item
One has $[D], [D'] \in \lambda^{\circ}$
for some GIT-chamber $\lambda \in \Lambda(\b{X},H)$.
\end{enumerate}
\end{proposition}

Let us have a look at the effect of blowing up
and more general modifications on the Cox ring.
In general, this is delicate, even finite 
generation may be lost.
We discuss a class of toric ambient 
modifications having good properties in this 
regard; we restrict to hypersurfaces,
for the general case see~\cite{Ha2}.

The starting point is a variety $X_0$ arising 
from a bunched ring $(R_0,\mathfrak{F}_0,\Phi_0)$
where the $K_0$-graded algebra $R_0$ is of the form 
$\KK[T_1, \ldots, T_r] / \bangle{f_0}$
and $\mathfrak{F}_0$ consists of 
the (pairwise nonassociated $K_0$-prime)
classes $T_i + \bangle{f_0}$; as usual, 
we denote their $K_0$-degrees by $w_i$.
We obtained the canonical toric embedding
via
$$ 
\xymatrix{
{\b{X}_0}
\ar@{}[r]|\supseteq
&
{\rq{X}_0}
\ar[r]
\ar[d]_{\quot H}
& 
{\rq{Z}_0}
\ar[d]^{\quot H}
\ar@{}[r]|\subseteq
&
{\b{Z}_0}
\\
&
X_0 
\ar[r]_{\imath}
&
Z_0
&
}
$$
The morphism $\rq{Z}_0 \to Z_0$ is given by a map 
$\rq{\Sigma}_0 \to \Sigma_0$ of fans
living in lattices $F_0 = \ZZ^r$ and $N_0$.
Let $v_1, \ldots, v_r$ be the primitive
lattice vectors in the rays of $\Sigma_0$
and suppose that for $2 \le d \le r$,
the cone $\sigma_0$ generated by 
$v_1, \ldots, v_d$ belongs to $\Sigma_0$.
Consider the stellar subdivision 
$\Sigma_1 \to \Sigma_0$ 
at a vector
\begin{eqnarray*}
v_\infty 
& = &
a_1v_1 + \cdots + a_dv_d.
\end{eqnarray*}
Let $m_\infty$ be the index of this subdivision,
i.e.~the gcd of the entries of $v_\infty$,
and denote the associated toric modification  by
$\pi \colon Z_1 \to Z_0$.
Then we obtain the strict transform 
$X_1 \subseteq Z_1$ mapping onto 
$X_0 \subseteq Z_0$.
Moreover, we have commutative diagrams
$$ 
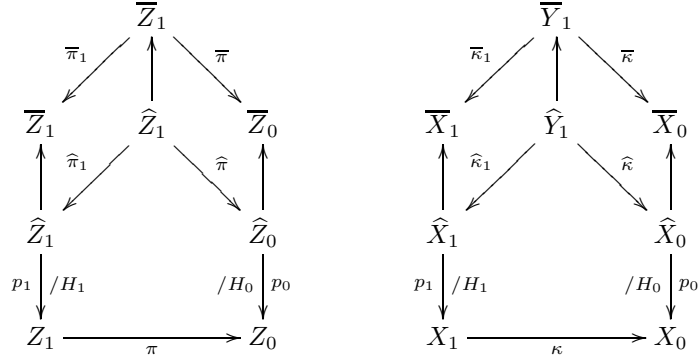
\xymatrix{
&
{\b{Z}_1}  \ar[dr]^{\b{\pi}} \ar[dl]_{\b{\pi}_1}
& 
\\
{\b{Z}_1}
&
{\rq{Z}_1} \ar[u] \ar[dl]_{\rq{\pi}_1} \ar[dr]^{\rq{\pi}}
& 
{\b{Z}_0}
\\
{\rq{Z}_1} \ar[d]_{p_1}^{/ H_1} \ar[u]
&
& 
{\rq{Z}_0} \ar[u] \ar[d]^{p_0}_{/ H_0}
\\
Z_1 
\ar[rr]_{\pi}
& &
Z_0
}
\qquad \qquad
\xymatrix{
&
{\b{Y}_1}  \ar[dr]^{\b{\kappa}} \ar[dl]_{\b{\kappa}_1}
& 
\\
{\b{X}_1}
&
{\rq{Y}_1} \ar[u] \ar[dl]_{\rq{\kappa}_1} \ar[dr]^{\rq{\kappa}}
& 
{\b{X}_0}
\\
{\rq{X}_1} \ar[d]_{p_1}^{/ H_1} \ar[u]
&
& 
{\rq{X}_0} \ar[u] \ar[d]^{p_0}_{/ H_0}
\\
X_1 
\ar[rr]_{\kappa}
& &
X_0
}
$$
where $H_1 = \Spec(\KK[K_1])$ for 
$K_1 = E_1 / M_1$ in analogy to 
$K_0 = E_0 / M_0$ etc..
Note that the Cox ring 
$\KK[T_1, \ldots, T_r,T_\infty]$
of $\b{Z}_1$ comes with a $K_1$-grading.
Moreover, with respect 
to the coordinates corresponding 
to the rays of the fans $\Sigma_i$, 
the map 
$\b{\pi} \colon \b{Z}_1 \to \b{Z}_0$ is given by 
$$ 
\b{\pi}(z_1, \ldots, z_r,z_\infty) 
\ = \
(z_\infty^{a_1}z_1, 
\ldots, z_\infty^{a_d}z_d, 
z_{d+1}, \ldots, z_r).
$$
We want to formulate an explicit condition
on the setting which guarantees that 
$\Gamma(\b{X}_1,\mathcal{O})$ is the 
Cox ring of the proper transform $X_1$.
For this, consider the grading 
$$ 
\KK[T_1, \ldots, T_r] 
\ = \ 
\bigoplus_{k \ge 0} \KK[T_1, \ldots, T_r]_k,
\quad
\text{where }
\deg(T_i) 
\ := \ 
\begin{cases}  
a_i & i \le d,
\\
0   & i \ge d+1.
\end{cases}
$$
Then we may write $f_0 = g_{k_0} + \cdots + g_{k_m}$
where $k_0 < \cdots < k_m$ and each $g_{k_i}$
is a nontrivial polynomial  having 
degree $k_i$ with respect to this grading.

\begin{definition}
We say that the polynomial $f_0 \in \KK[T_1, \ldots, T_r]$ 
is {\em admissible\/} if 
\begin{enumerate}
\item
the toric orbit $0 \times \TT^{r-d}$ intersects 
$\b{X}_0 = V(f_0)$,
\item
$g_{k_0}$ 
is a $K_1$-prime polynomial in at 
least two variables.
\end{enumerate}
\end{definition}

Note that for the case of a free abelian group $K_1$,
the second condition just means that $g_{k_0}$ is 
an irreducible polynomial.

\begin{proposition}
\label{singlerel}
If, in the above setting, the 
polynomial $f_0$ is admissible,
then the Cox ring of the strict transform 
$X_1 \subseteq Z_1$ is
$$
\mathcal{R}(X_1) 
\ = \ 
\KK[T_1, \ldots, T_r,T_\infty] / 
\bangle{f_1(T_1, \ldots, T_r,\sqrt[m_\infty]{T_\infty})},
$$
where in
$$ 
f_1 
\ := \ 
\frac{f_0(T_\infty^{a_1}T_1, 
\ldots, T_\infty^{a_d}T_d, 
T_{d+1}, \ldots, T_r )}{T_\infty^{k_0}}
\ \in \ 
\KK[T_1, \ldots, T_r,T_\infty]
$$
only powers $T_\infty^{lm_\infty}$ with $l \ge 0$
of $T_\infty$ occur, 
and the notation $\sqrt[m_\infty]{T_\infty}$ means 
replacing 
 $T_\infty^{lm_\infty}$ with $T_\infty^l$ in $f_1$. 
\end{proposition}

\begin{example}
Consider the polynomial 
$f_0 := T_1 + T_2^2+T_3T_4$ and the 
factorial algebra 
$R_0 := \KK[T_1,\ldots,T_4]/f_0$.
Then $R_0$ is graded by $K_0 := \ZZ$ 
via the weight matrix 
\begin{eqnarray*}
Q_0 
& = & 
[2,1,1,1].
\end{eqnarray*}
With $\mathfrak{F}_0 := (\b{T}_1, \ldots, \b{T}_4)$ 
and $\Phi_0 = \{\QQ_{\ge 0}\}$, 
we obtain a bunched ring $(R_0,\mathfrak{F}_0,\Phi_0)$.
The associated variety $X_0$ is a projective surface.
In fact, there is an isomorphism $X_0 \to \PP_2$ 
induced by
$$ 
\KK^4 \ \to \ \KK^3,
\qquad
(z_1,z_2,z_3,z_4)
\ \mapsto \ 
(z_2,z_3,z_4).
$$
The canonical toric ambient variety $Z_0$ of $X_0$
is an open subset of the weighted projective space
$\PP(2,1,1,1)$.
To obtain a defining fan $\Sigma_0$ of $Z_0$, 
consider the Gale dual matrix
\begin{eqnarray*}
P_0 
& := & 
\left(
\begin{array}{rrrr}
-1 & 2 & 0 & 0
\\
-1 & 0 & 1 & 1
\\
0 & 1 & -1 & 0 
\end{array}
\right).
\end{eqnarray*}
Let $v_1,\ldots,v_4$ denote its columns $v_1,\ldots,v_4$. 
According to~\ref{ex:bunchedpolring}, the maximal cones of 
$\Sigma_0$ are $P(\gamma_0^\perp \cap \delta)$,
where $\gamma_0$ runs through the covering collection
$\cov(\Phi_0)$.
In terms of the columns $v_1,\ldots,v_4$ of $P$ they 
are given as  
$$ 
\sigma_{1,2,3} \ := \ \cone(v_1,v_2,v_3),
\qquad
\sigma_{1,2,4} \ := \ \cone(v_1,v_2,v_4),
\qquad
\sigma_{3,4} \ := \ \cone(e_3,e_4).
$$
Now we subdivide $\sigma_{1,2,3}$ by inserting 
the ray through $v_\infty = 3v_1+v_2+2v_3$.
Note that $f_0$ is admissible; the $g_{k_0}$-term is 
$T_2^2+T_3T_4$.
According to~\ref{singlerel}, the defining equation
of the Cox ring of the proper transform $X_1$ is
$$
f_1
\ = \ 
\frac{T_\infty^3T_1+T_\infty^2T_2^2+T_\infty^2T_3T_4}{T_\infty^2}
\ = \ 
T_\infty T_1 + T_2^2 + T_3T_4.
$$ 
In order to obtain the grading matrix, we have to 
look at the matrix with the new primitive generators
as columns. In abuse of notation, we put  
$v_{\infty} = (-1,-1,-1)$ at the first place:
\begin{eqnarray*}
P_1 
& := & 
\left(
\begin{array}{rrrrr}
-1 & -1 & 2 & 0 & 0
\\
-1 & -1 & 0 & -1 & 0
\\
-1 & 0 & 1 & -1 & 0 
\end{array}
\right).
\end{eqnarray*}
The degree matrix is the Gale dual: 
\begin{eqnarray*}
Q_1 
& = & 
\left(
\begin{array}{rrrrr}
1 & -1 & 0 & -1 & 1
\\
1 & 1 & 1 & 0 & 2
\end{array}
\right).
\end{eqnarray*}
Thus, up to renaming of the variables, we obtain 
the Cox ring of the singular del Pezzo surface 
considered in~\ref{ex:delpezzo} and~\ref{ex:delpezzo2}.
In particular, we see that this del Pezzo surface 
is a modification of the projective plane. 
\end{example}

\section{Third lecture}

\subsection{Varieties with torus action}
We describe the Cox ring of a variety with 
torus action following~\cite{HaSu}.
First recall the example of a complete 
toric variety $X$.
Its Cox ring is given in terms of the 
prime divisors $D_1, \ldots, D_r$ in 
the boundary $X \setminus T \mal x_0$ of the 
open orbit:
$$
\mathcal{R}(X)
\ = \ 
\KK[T^{D_1}, \ldots, T^{D_r}],
\qquad
\deg(T^{D_i}) \ = \ [D] \ \in \ \Cl(X).
$$

Now let $X$ be any normal complete 
variety with finitely generated divisor class
group and consider an effective algebraic torus 
action $T \times X \to X$, where $\dim(T)$
may be less than $\dim(X)$.
For a point $x \in X$, denote by 
$T_x \subseteq T$ its isotropy group.
The points with finite isotropy group form 
a non-empty $T$-invariant open subset
$$ 
X_0
\ := \ 
\{x \in X; \; T_x \text{ is finite} \}
\ \subseteq \ 
X.
$$
This set will replace the open orbit of a 
toric variety.
Let $E_k$, where $1 \le k \le m$, denote 
the prime divisors in $X \setminus X_0$;
note that each $E_k$ is $T$-invariant with
infinite generic isotropy, i.e.~the subgoup 
of $T$ acting trivially on $E_k$ is infinite.
According to a result of Sumihiro~\cite{Su}, 
there is a geometric quotient 
$q \colon X_0 \to X_0 / T$ 
with an irreducible normal but possibly 
non-separated orbit space 
$X_0 / T$.

\begin{example}
\label{ex:blp1p1}
Consider $\PP_1 \times \PP_1$ the 
$\KK^*$-action given w.r.t. inhomogeneous 
coordinates by $t \mal (z,w)  = (z,tw)$.
Let $X$ be the $\KK^*$-equivariant blow up
of $\PP_1 \times \PP_1$ at the fixed points 
$(0,0)$, $(1,0)$ and $(\infty,0)$.
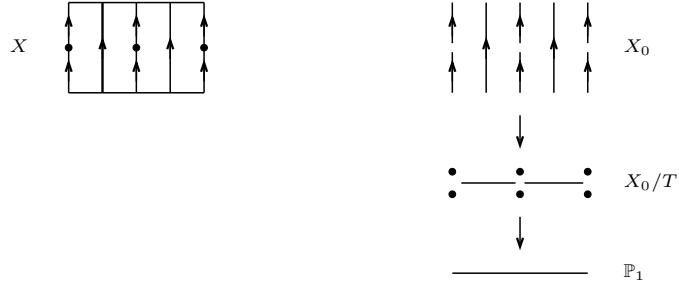
\begin{center}

\medskip

\begin{picture}(0,0)%
\includegraphics{blp2p2cox.pstex}%
\end{picture}%
\setlength{\unitlength}{1243sp}%
\begingroup\makeatletter\ifx\SetFigFontNFSS\undefined%
\gdef\SetFigFontNFSS#1#2#3#4#5{%
  \reset@font\fontsize{#1}{#2pt}%
  \fontfamily{#3}\fontseries{#4}\fontshape{#5}%
  \selectfont}%
\fi\endgroup%
\begin{picture}(12270,5632)(-284,-5210)
\put(-269,-601){\makebox(0,0)[lb]{\smash{{\SetFigFontNFSS{7}{8.4}{\familydefault}{\mddefault}{\updefault}{\color[rgb]{0,0,0}$X$}%
}}}}
\put(11971,-601){\makebox(0,0)[lb]{\smash{{\SetFigFontNFSS{7}{8.4}{\familydefault}{\mddefault}{\updefault}{\color[rgb]{0,0,0}$X_0$}%
}}}}
\put(11971,-3301){\makebox(0,0)[lb]{\smash{{\SetFigFontNFSS{7}{8.4}{\familydefault}{\mddefault}{\updefault}{\color[rgb]{0,0,0}$X_0/T$}%
}}}}
\put(11971,-5101){\makebox(0,0)[lb]{\smash{{\SetFigFontNFSS{7}{8.4}{\familydefault}{\mddefault}{\updefault}{\color[rgb]{0,0,0}$\PP_1$}%
}}}}
\end{picture}%

\medskip

\end{center}
The open set $X_0$ is obtained by removing the 
two fixed point curves and the three isolated fixed 
points.
The quotient $X_0/T$ is the non-separated projective 
line with the points $0,1,\infty$ doubled; note that 
there is a canonical map $X_0/T \to \PP_1$.
\end{example}

As it turns out, one always finds 
kind of {\em separation\/} for the orbit space
in our setting: 
there are a rational map 
$\pi \colon X_0/T \dashrightarrow Y$,
an open subset $W \subseteq X_0/T$ 
and prime divisors $C_0, \ldots, C_r$ on $Y$ 
such that following holds:
\begin{itemize}
\item 
the complement of $W$ in $X_0/T$ is of codimension
at least two and the restriction 
$\pi \colon W \to X_0/T$ is a local isomorphism,
\item
each inverse image $\pi^{-1}(C_i)$ is a disjoint 
union of prime divisors $C_{ij}$, where $1 \le j \le n_i$, 
\item 
the map $\pi$ is an isomorphism over 
$Y \setminus (C_0 \cup \ldots \cup C_r)$
and all prime divisors of $X_0$ 
with nontrivial generic isotropy occur among 
the $D_{ij} := q^{-1}(C_{ij})$.
\end{itemize}
Let $l_{ij} \in \ZZ_{\ge 1}$ be the order
of the generic isotropy group of the $T$-action 
on the prime divisor $D_{ij}$
and define monomials 
$T_i^{l_i} := T_{i1}^{l_{i1}} \cdots T_{in_i}^{l_{in_i}}$
in the variables $T_{ij}$.
Moreover, let $1_{E_k}$ and $1_{D_{ij}}$ denote the 
canonical sections of $E_k$ and $D_{ij}$.

\begin{theorem}
\label{fingenchar2}
There is a graded injection 
$\mathcal{R}(Y) \to \mathcal{R}(X)$ 
of Cox rings and $S_k \mapsto 1_{E_k}$,
$T_{ij} \mapsto 1_{D_{ij}}$ defines
an isomorphism of $\Cl(X)$-graded rings
\begin{eqnarray*}
\mathcal{R}(X)
& \cong &
\mathcal{R}(Y)[T_{ij},S_1,\ldots,S_m; \; 
0 \le i \le r, \, 1 \le j \le n_i]
\ / \ 
\bangle{T_i^{l_i} - 1_{C_i}; \; 0 \le i \le r},
\end{eqnarray*}
where the $\Cl(X)$-grading on the right hand side 
is defined by associating to $S_k$ the class of $E_k$ 
and to $T_{ij}$ the class of $D_{ij}$.
\end{theorem}

As a direct consequence, we obtain that finite 
generation of the Cox ring of $X$ is determined 
by the separation $Y$ of the orbit space $X_0/T$.

\begin{corollary}
The Cox ring $\mathcal{R}(X)$ is finitely
generated if and only if $\mathcal{R}(Y)$ 
is so. 
\end{corollary}

We specialize to the case that the $T$-action 
on $X$ is of {\em complexity\/} one, i.e.~its biggest 
orbits are of codimension one in $X$.
Then the orbit space $X_0/T$ is of dimension
one and smooth.

\begin{remark}
For a normal complete variety $X$ with  
a torus action $T \times X \to X$ 
of complexity one, the following
statements are equivalent.
\begin{enumerate}
\item
$\Cl(X)$ is finitely generated.
\item
$X$ is rational.
\end{enumerate}
Moreover, if one of these statements holds, then 
the separation of the orbit space is a morphism 
$\pi \colon X_0/T \to \PP_1$.
\end{remark}

The former prime divisors $C_i \subseteq Y$ 
are now points $a_0, \ldots, a_r \in \PP_1$ and
$\pi^{-1}(a_i)$ consists of points 
$x_{i1}, \ldots, x_{in_i}  \in X_0/T$.
As before, we may assume that all prime divisors 
of $X_0$ with 
non trivial generic isotropy occur among 
the prime divisors $D_{ij} = q^{-1}(x_{ij})$.
Consider again the monomials 
$T_i^{l_i} = T_{i1}^{l_{i1}} \cdots T_{in_i}^{l_{in_i}}$,
where $l_{ij} \in \ZZ_{\ge 1}$ is 
the order of the generic isotropy group of 
$D_{ij}$.
Write $a_i = [b_i,c_i]$ with 
$b_i,c_i \in \KK$.
For every $0 \le i \le r-2$ set 
$j := i+1$, $k := i+2$ and define a trinomial
\begin{eqnarray*}
g_i 
& := &
(b_jc_k - c_jb_k)T_i^{l_i}
\ + \ 
(b_kc_i - c_kb_i)T_j^{l_j}
\ + \ 
(b_ic_j - c_ib_j)T_k^{l_k}.
\end{eqnarray*}

\begin{theorem}
\label{complexity1}
Let $X$ be a normal complete variety 
with 
finitely generated divisor class group
and an effective algebraic torus action 
$T \times X \to X$ 
of complexity one.
Then, in terms of the data defined
above, the Cox ring of $X$ is given as
\begin{eqnarray*}
\mathcal{R}(X)
&  \cong &
\KK[S_1,\ldots, S_m, T_{ij}; \; 0 \le i \le r, \; 1 \le j \le n_i] 
\ / \
\bangle{g_i; \; 0 \le i \le r-2}.
\end{eqnarray*}
where $1_{E_k}$ corresponds to $S_k$,
and $1_{D_{ij}}$ to $T_{ij}$,
and the $\Cl(X)$-grading on the right hand side 
is defined by associating to $S_k$ the class 
of $E_k$ and to $T_{ij}$ the class of $D_{ij}$.
\end{theorem}

Let us discuss some applications of 
Theorem~\ref{complexity1}.
In a first application we compute 
the Cox ring of a surface 
obtained by repeated blowing up 
points of the projective plane 
that lie on a given line; the case
$n_0 = \ldots = n_r = 1$ was 
done by other methods in~\cite{Ot}.

\begin{example}[Blowing up points on a line]
Consider a line $Y \subseteq \PP_2$
and points $p_0, \ldots, p_r \in Y$.
Let $X$ be the surface obtained by 
blowing up $n_i$ times the point $p_i$,
where $0 \le i \le r$; 
in every step, we identify $Y$ with its proper transform 
and $p_i$ with the point in the intersection 
of $Y$ with the exceptional curve. Set
\begin{eqnarray*}
g_i 
& := & 
(b_jc_k - c_jb_k)T_i
\ + \ 
(b_kc_i - c_kb_i)T_j
\ + \ 
(b_ic_j - c_ib_j)T_k,
\end{eqnarray*}
with $p_i = [b_i,c_i]$ in $Y = \PP_1$,
the monomials $T_i = T_{i0} \cdots T_{in_i}$
and the indices 
$k = i+2$,  $j= i+1$.
Then the Cox ring of the surface $X$ is given as
\begin{eqnarray*}
\mathcal{R}(X) 
& = & 
\KK[T_{ij},S; \; 0 \le i \le r, \; 0 \le j \le n_i]
\ / \ 
\bangle{g_i; \; 0 \le i \le r-2}.
\end{eqnarray*}
We verify this using a $\KK^*$-action; note that
blowing up $\KK^*$-fixed points always can be 
made equivariant.
With respect to suitable homogeneous coordinates
$z_0,z_1,z_2$, we have $Y = V(z_0)$.
Let $\KK^*$ act via 
\begin{eqnarray*}
t \mal [z] 
& := & 
[z_0,tz_1,tz_2].
\end{eqnarray*}
Then $E_1 := Y$ is a fixed point curve, 
the $j$-th (equivariant) blowing up of $p_i$ 
produces an invariant exceptional divisor 
$D_{ij}$ with a free $\KK^*$-orbit inside
and Theorem~\ref{complexity1} gives the 
claim.
\end{example}

As a further application we indicate 
a general recipe for computing the Cox ring 
of a rational hypersurface given by a 
trinomial equation in the projective space.
We perform this for the $E_6$ cubic 
surface in $\PP_3$; the Cox ring 
of the resolution of this surface 
has been computed in~\cite{HaTs}.

\begin{example}[The $E_6$ cubic surface]
\label{ex:E6cubic}
There is a cubic surface~$X$ in the projective 
space having singular locus $X^{\rm sing} = \{x_0\}$ 
and $x_0$ of type $E_6$.
The surface is unique up to projectivity and 
can be realized as follows:
$$  
X 
\ = \ 
V(z_1z_2^2 + z_2z_0^2 + z_3^3)
\ \subseteq \ 
\PP_3.
$$
Note that the defining equation is a trinomial 
but not of the shape of those occurring in 
Theorem~\ref{complexity1}.
However, any trinomial hypersurface in a projective
space comes with a complexity one torus action.
Here, we have the $\KK^*$-action
\begin{eqnarray*}
t \mal [z_0,\ldots,z_4] 
& = & 
[z_0,t^{-3}z_1,t^3z_2,tz_3].
\end{eqnarray*}
This allows  us to use Theorem~\ref{complexity1} for
computing the Cox ring.
The task is to find the divisors $E_k,D_{ij}$ and the 
orders $l_{ij}$ of the isotropy groups.
Note that $\KK^*$ acts freely on the big torus of
$\PP_3$.
The intersections of~$X$ with the toric prime 
divisors $V(z_i) \subseteq \PP_3$ are given as
$$ 
X \cap V(z_0)
\ = \ 
V(z_0, z_1z_2^2+z_3^3),
\qquad
X \cap V(z_1)
\ = \ 
V(z_1, z_2z_0^2+z_3^3),
$$
$$ 
X \cap V(z_3)
\ = \ 
V(z_3,z_2(z_1z_2+z_0^2))
\ = \ 
(X \cap V(z_2)) \cup (X \cap V(z_1z_2+z_0^2)).
$$
The first two sets are irreducible and 
both of them intersect the big torus orbit
of the respective toric prime divisors
$V(z_0)$ and $V(z_1)$.
In order to achieve this also  
for $V(z_2)$, $V(z_3)$, 
we use a suitable weighted blow 
up of $\PP_3$ at $V(z_2) \cap V(z_3)$.
In terms of fans, this means to perform 
a certain stellar subdivision.
Consider the matrices
$$ 
P
\ = \
\left(
\begin{array}{rrrr}
-1 & 1 & 0 & 0 
\\
-1 & 0 & 1 & 0
\\
-1 & 0 & 0 & 1
\end{array}
\right),
\qquad
\qquad
P'
\ = \
\left(
\begin{array}{rrrrr}
-1 & 1 & 0 & 0 & 0 
\\
-1 & 0 & 1 & 0 & 3
\\
-1 & 0 & 0 & 1 & 1
\end{array}
\right).
$$
The columns $v_0, \ldots, v_3$ of $P$ 
are the primitive generators of the fan 
$\Sigma$ of $Z := \PP_3$ and we obtain a 
fan $\Sigma'$ subdividing $\Sigma$ at 
the last column  $v_4$ of $P'$;
note that $v_4$ is located on the 
tropical variety ${\rm trop}(X)$.
Consider the associated toric morphism
and the proper transform
$$ 
\pi \colon Z \ \to \ Z,
\qquad \qquad 
X' \ := \  \b{\pi^{-1}(X \cap \TT^3)} \ \subseteq \ Z_1.
$$
A simple computation shows that 
the intersection of $X'$ with the 
toric prime divisors of $Z'$ is 
irreducible and intersects their 
big orbits.
Moreover, $\pi \colon X' \to X$ is an
isomorphism, because along $X'$ nothing
gets contracted.
To proceed, note that we have no divisors 
of type $E_k$ and that there is a
commutative diagram
$$ 
\xymatrix{
X'_0 
\ar[r]
\ar[d]_{/ \KK^*}
&
Z'_0
\ar[d]^{/ \KK^*}
\\
X'_0 / \KK^*
\ar[r]
&
Z'_0 / \KK^*
}
$$
We determine the quotient $Z'_0 \to Z'_0 / \KK^*$.
The group $\KK^*$ acts on $Z'$ via homomorphism
$\lambda_v \colon \KK^* \to \TT^3$ corresponding 
to $v = (-3,3,1) \in \ZZ^3$.
The quotient by this action is the toric morphism
given by any map $S \colon \ZZ^3 \to \ZZ^2$ 
having $\ZZ \mal v$ as its kernel.
We take $S$ as follows and compute the images 
of the columns $v_0,\ldots,v_4$ of $P'$:
$$ 
\qquad
S 
\ := \  
\left(
\begin{array}{rrr}
1 & 0 & 3
\\
0 & 1 & -3
\end{array}
\right),
\qquad \qquad
S \cdot P' 
\ = \  
\left(
\begin{array}{rrrrr}
-4 & 1 & 0 & 3 & 3 
\\
2 & 0 & 1 & -3 & 0
\end{array}
\right).
$$
This shows that the toric divisors $D_Z^1$ and $D^4_Z$ 
corresponding to $v_1$ and $v_4$ are mapped to a 
doubled divisor in the non-separated quotient;
see also~\cite{ACHa}.
The generic isotropy group
of the $\KK^*$-action along the toric divisor
$D_Z^i$ is given as the 
gcd $l_i$ of the entries of the $i$-th column
of $S \mal P'$.
We obtain
$$ 
l_0 \ = \ 2,
\qquad
l_1 \ = \ 1,
\qquad
l_2 \ = \ 1,
\qquad
l_3 \ = \ 3,
\qquad
l_4 \ = \ 3.
\qquad
$$
By construction, the divisors $D_X^i := D_Z^i$
of the embedded variety $X'_0 \subseteq Z'_0$ 
inherit the orders $l_i$ of the isotropy groups 
and the behaviour with respect to the quotient map 
$X'_0 \to X'_0 / \KK^*$. Renaming these divisors 
by
$$
D_{01} := D_X^1, 
\quad
D_{02} := D_X^4,
\quad
D_{11} := D_X^3,
\quad
D_{21} := D_X^0,
$$ 
we arrive in the setting of Theorem~\ref{complexity1};
since $D_X^2$ has trivial isotropy group and 
gets not multiplied by the quotient map, it 
does not occur here.
The Cox ring of $X \cong X'$ is then given as 
\begin{eqnarray*}
\mathcal{R}(X)
& \cong & 
\KK[T_{01},T_{02},T_{11},T_{21}] / \bangle{T_{01}T_{02}^3 + T_{11}^3 + T_{21}^2}.
\end{eqnarray*}
\end{example}

\subsection{Writing down all Cox rings}
As we observed, the Cox ring of a rational 
complete normal variety with a complexity one 
torus action admits a nice presentation by 
trinomial relations.
Now we ask which of these trinomial rings 
occur as a Cox ring. 
First, we formulate the answer in algebraic
terms and then turn to the geometric point 
of view; for the details, we refer to~\cite{HaHe}.

\goodbreak

\begin{construction}
\label{constr:RAP}
Fix $r \in \ZZ_{\ge 1}$, a sequence 
$n_0, \ldots, n_r \in \ZZ_{\ge 1}$, set 
$n := n_0 + \ldots + n_r$, and fix  
integers $m \in \ZZ_{\ge 0}$ and $0 < s < n+m-r$.
The input data are 
\begin{itemize}
\item 
a {\em sequence\/} $A = (a_0, \ldots, a_r)$ 
of vectors $a_i = (b_i,c_i)$ in $\KK^2$
such that any pair $(a_i,a_j)$ with
$j \ne i$ is linearly independent,
\item 
an {\em integral block matrix\/} $P$ of size 
$(r + s) \times (n + m)$ the columns 
of which are pairwise different primitive
vectors generating $\QQ^{r+s}$ as a cone: 
\begin{eqnarray*}
P
& = & 
\left( 
\begin{array}{cc}
P_0 & 0 
\\
d & d'  
\end{array}
\right),
\end{eqnarray*}
where $d$ is an $(s \times n)$-matrix, $d'$ an $(s \times m)$-matrix 
and $P_0$ an $(r \times n)$-matrix build from
tuples $l_i := (l_{i1}, \ldots, l_{in_i}) \in \ZZ_{\ge 1}^{n_i}$ 
in the following way 
\begin{eqnarray*}
P_0
& = & 
\left(
\begin{array}{cccc}
-l_0 & l_1 &   \ldots & 0 
\\
\vdots & \vdots   & \ddots & \vdots
\\
-l_0 & 0 &\ldots  & l_{r} 
\end{array}
\right).
\end{eqnarray*}
\end{itemize}
Now we associate to any such pair $(A,P)$ a ring,
graded by $K := \ZZ^{n+m}/\rm{im}(P^*)$, where 
$P^*$ is the transpose of $P$.
For every $0 \le i \le r$, define a monomial
\begin{eqnarray*}
T_i^{l_i} 
&  := &
T_{i1}^{l_{i1}} \cdots T_{in_i}^{l_{in_i}}.
\end{eqnarray*}
Moreover, for any two indices $0 \le i,j \le r$,
set $\alpha_{ij} :=  \det(a_i,a_j)  =  b_ic_j-b_jc_i$
and for any three indices 
$0 \le i < j < k \le r$ define 
a trinomial
\begin{eqnarray*}
g_{i,j,k} 
& := &
\alpha_{jk}T_i^{l_i} 
\ + \ 
\alpha_{ki}T_j^{l_j} 
\ + \ 
\alpha_{ij}T_k^{l_k}.
\end{eqnarray*}
Note that all trinomials $g_{i,j,k}$ are $K$-homogeneous
of the same degree.
Setting $g_i := g_{i,i+1,i+2}$, we obtain a $K$-graded 
factor algebra 
\begin{eqnarray*}
R(A,P)
& := &
\KK[T_{ij},S_k; \; 0 \le i \le r, \, 1 \le j \le n_i, 1 \le k \le m] 
\ / \
\bangle{g_i; \; 0 \le i \le r-2}.
\end{eqnarray*}
\end{construction}

\begin{remark}
The polynomials $g_{i,j,k}$ can be written as determinants
in the following way:
\begin{eqnarray*}
g_{i,j,k}
& = & 
\det
\left(
\begin{array}{ccc}
b_i & b_j & b_k
\\
c_i & c_j & c_k 
\\
T_i^{l_i} & T_j^{l_j} & T_k^{l_k}
\end{array}
\right).
\end{eqnarray*}
\end{remark}

\begin{theorem}
\label{thm:RAPalgebras}
Let $(A,P)$ be data as in~\ref{constr:RAP}.
Then the algebra $R := R(A,P)$ is a normal 
factorially $K$-graded complete intersection; 
we have $R^* = \KK^*$, the $K$-grading is 
almost free and $R_0 = \KK$ holds.
Moreover, the variables $T_{ij}$, $S_k$ define 
a system of pairwise nonassociated $K$-prime 
generators for $R$.
\end{theorem}

We say that the pair $(A,P)$ 
is {\em sincere\/}, if $r \ge 2$ and 
$n_il_{ij} > 1$ for all $i,j$ hold;
this ensures that there exist in fact relations 
$g_{i,j,k}$ and none of these relations 
contains a linear term.
The following statement tells us which of the 
$R(A,P)$ are factorial.

\begin{theorem}
\label{thm:RAPUFD}
Let $(A,P)$ be a sincere pair. Then the following 
statements are equivalent.
\begin{enumerate}
\item
The algebra $R(A,P)$ is a unique factorization 
domain.
\item
The group $\ZZ^r/ {\rm im}(P_0)$ is torsion free.
\item
The numbers $\gcd(l_i)$ and $\gcd(l_j)$ are 
coprime for any $0 \le i < j \le r$.
\end{enumerate}
\end{theorem}

More generally, one can in a similar manner 
to~\ref{thm:RAPalgebras} determine all affine 
algebras $R$ with an effective factorial $K$-grading 
of complexity one such that $R_0 = \KK$ holds, 
see~\cite{HaHe}.
The characterization of the factorial ones 
is the same as in~\ref{thm:RAPUFD}.
In dimensions two and three, the factorial
algebras with a complexity one multigrading
are described in early work of Mori~\cite{Mo}
and Ishida~\cite{Is}.

\begin{example}
The algebra 
$\KK[T_{01},T_{11},T_{21}] / \bangle{T_{01}^2+T_{11}^2+T_{21}^2}$
becomes graded by the group
$K = \ZZ \oplus  \ZZ/2\ZZ \oplus  \ZZ/2\ZZ$
via
$$ 
\deg(T_{01}) = (1,\b{0}, \b{0}),
\qquad
\deg(T_{11}) = (1,\b{1}, \b{0}),
\qquad
\deg(T_{21}) = (1,\b{0}, \b{1}).
$$
This is an effective factorial grading 
of complexity one.
However, the grading is not almost free.
Thus, the algebra is not a Cox ring.
\end{example}

Let us turn to geometric aspects.
We want to see the effective complexity 
one torus action on the varieties 
having a Cox ring $R(A,P)$.
Existence of this action can be obtained 
by looking at the maximal possible 
grading of $R(A,P)$. 
We use here the canonical toric embedding 
which provides a little more geometric
information.

\begin{construction}
\label{constr:RAP2Taction}
Take a $K$-graded algebra $R = R(A,P)$ 
as constructed in~\ref{constr:RAP}
and let $\mathfrak{F}$ be the system of 
pairwise nonassociated $K$-prime 
generators defined by the variables 
$T_{ij}$ and $S_k$.
Given any $\mathfrak{F}$-bunch $\Phi$, 
we obtain a bunched ring $(R,\mathfrak{F},\Phi)$
and an associated variety $X$.
Consider the mutually dual sequences
$$ 
\xymatrix@R=10pt{
0 \ar[r]
&
L
\ar[r]
&
F
\ar[r]^{P}
&
N,
&
\\
0 \ar@{<-}[r]
&
K
\ar@{<-}[r]_{Q}
&
E
\ar@{<-}[r]_{P^*}
&
M
\ar@{<-}[r]
&
0.
}
$$
Recall that $Q \colon E \to K$ defines 
the $K$-degrees of the variables $T_{ij}$, $S_k$ 
and note that $P$ is indeed our defining 
matrix of the algebra $R(A,P)$.
Now, let $\rq{\Sigma}$ denote the fan in $F$ 
generated by $\Phi$ and 
let $\Sigma$ be its quotient fan in $N$;
its rays have the columns $v_{ij}$ and $v_k$ 
of $P$ as their primitive generators
Moreover, consider
$$ 
\b{X} 
\ := \
V(g_0, \ldots, g_{r-2})
\ \subseteq \ 
\KK^{n+m}.
$$
The fan $\rq{\Sigma}$ defines an open toric subvariety
$\rq{Z} \subseteq \KK^{n+m}$ and the toric morphism 
$p \colon \rq{Z} \to Z$ defined by $P \colon F \to N$ 
onto the toric variety $Z$ associated to $\Sigma$ 
is a characteristic space. 
The canonical toric embedding of $X$ was obtained
via the commutative diagram
$$ 
\xymatrix{
{\b{X}}
\ar@{}[r]|\supseteq
&
{\rq{X}}
\ar[r]
\ar[d]_{\quot H}
&
{\rq{Z}}
\ar@{}[r]|\subseteq
\ar[d]^{\quot H}
&
{\KK^{n+m}}
\\
& 
X
\ar[r]
& 
Z
& 
}
$$
Let $T_Z := \Spec \, \KK[M]$ be the acting torus of 
$Z$ and let $T \subseteq T_Z$ be the subtorus 
corresponding to the inclusion $0 \times \ZZ^s \to N$.
Then $T$ acts on $Z$ leaving $X \subseteq Z$ 
invariant and $T \times X \to X$ is an effective 
complexity one action.
We also write $X = X(A,P,\Phi)$ for this $T$-variety.
\end{construction}

\begin{theorem}
Let $X$ be a normal complete $A_2$-variety with 
an effective complexity one torus action.
Then $X$ is equivariantly isomorphic
to a variety $X(A,P,\Phi)$
constructed in~\ref {constr:RAP2Taction}
\end{theorem}

In particular, this enables us to apply the 
language of bunched rings to varieties with 
a complexity one torus action. 
The following is an application of~\ref{thm:ample} 
and~\ref{thm:canondiv}.

\begin{proposition}
\label{cor:canondiv}
Let $X$ be a complete normal rational variety with an 
effective algebraic torus action $T \times X \to X$
of complexity one.
\begin{enumerate}
\item
The cone of divisor classes without fixed components
is given by
$$
\qquad \qquad
\bigcap_{1 \le k \le m}
\cone([E_s], [D_{ij}]; \; s \ne k) 
\ \cap \ 
\bigcap_{\genfrac{}{}{0pt}{}{0 \le i \le r}{1 \le j \le n_i}}
\cone([E_k], [D_{st}]; \; (s,t) \ne (i,j)).
$$
\item 
For any $0 \le i \le r$, one obtains a canonical 
divisor for $X$ by 
$$
\max(0,r-1) \cdot \sum_{j=0}^{n_i} l_{ij} D_{ij} 
\ - \ 
\sum_{k=1}^m E_k
\ -  \ 
\sum_{i,j} D_{ij}.
$$
\end{enumerate}
\end{proposition}

The explicit description of rational varieties 
with a complexity one torus action by the data 
$(A,P)$  can be used for classifications.
We will present a result of~\cite{HaHeSu}
on Fano threefolds with free class group 
of rank one; we first note a general observation
on varieties of this type made there.

Let $X$ be an arbitrary 
complete $d$-dimensional
variety with divisor class group
$\Cl(X) \cong \ZZ$; thus, we do not 
require the presence of a torus action.
The Cox ring $\mathcal{R}(X)$
is finitely generated and the total
coordinate space $\b{X} := \Spec \, \mathcal{R}(X)$
is a factorial affine variety coming
with an action of $\KK^*$ defined by
the $\Cl(X)$-grading of $\mathcal{R}(X)$.
Choose a system $f_1, \ldots, f_\nu$ of
homogeneous pairwise nonassociated
prime generators for $\mathcal{R}(X)$.
This provides a $\KK^*$-equivariant
embedding
$$
\b{X} \ \to \ \KK^{\nu},
\qquad
\b{x} \ \mapsto \ (f_1(\b{x}), \ldots, f_{\nu}(\b{x})),
$$
where $\KK^*$ acts  diagonally with the
weights $w_i = \deg(f_i) \in \Cl(X) \cong \ZZ$
on $\KK^{\nu}$.
Moreover, $X$ is the geometric
$\KK^*$-quotient of
$\rq{X} := \b{X} \setminus \{0\}$,
and the quotient map
$p \colon \rq{X} \to X$ is a characteristic 
space.

\begin{proposition}
\label{Prop:FanoPicard}
For any
$\b{x} =(\b{x}_1,\ldots,\b{x}_{\nu}) \in \rq{X}$
the local divisor class group $\Cl(X,x)$
of $x := p(\b{x})$
is finite of order $\gcd(w_i; \; \b{x}_i \ne 0)$.
The index of the Picard group $\Pic(X)$ in
$\Cl(X)$ is given by
\begin{eqnarray*}
[\Cl(X):\Pic(X)]
& = &
\mathrm{lcm}_{x \in X}( |\Cl(X,x)| ).
\end{eqnarray*}
Suppose that the ideal of $\b{X} \subseteq \KK^{\nu}$
is generated by $\Cl(X)$-homogeneous
polynomials $g_1, \ldots, g_{\nu-d-1}$
of degree $\gamma_j := \deg(g_j)$.
Then one obtains
$$
-\mathcal{K}_X
\ = \
\sum_{i=1}^{\nu} w_i - \sum_{j=1}^{\nu-d-1} \gamma_j,
\qquad
(-\mathcal{K}_X )^d
\ = \
\left(\sum_{i=1}^{\nu} w_i - \sum_{j=1}^{\nu-d-1} \gamma_j\right)^d
\frac{\gamma_1 \cdots \gamma_{\nu-d-1}}{w_1 \cdots w_\nu}
$$
for the anticanonical class
$-\mathcal{K}_X \in \Cl(X) \cong \ZZ$.
In particular, $X$ is a Fano variety
if and only if the following inequality holds
\begin{eqnarray*}
 \sum_{j=1}^{\nu-d-1} \gamma_j
 & < &
 \sum_{i=1}^{\nu} w_i.
\end{eqnarray*}
\end{proposition}

Combining this with the explicit description in the 
presence of a complexity one torus action, 
one obtains bounds for the defining data.
This leads to classification results.
We present here the results for locally factorial 
threefolds; see~\cite{HaHeSu} for details
and more.

\begin{theorem}
The following table lists the 
Cox rings $\mathcal{R}(X)$ 
of the three-dimensio\-nal 
locally factorial non-toric Fano 
varieties $X$ with an effective two torus 
action and $\Cl(X) = \ZZ$.
\begin{center}
\begin{longtable}[htbp]{llll}
\toprule
No.
&
$\mathcal{R}(X)$ 
& 
$(w_1,\ldots, w_5)$
& 
$(-K_X)^3$
\\
\midrule
1 
\hspace{.5cm}
&
$
\KK[T_1, \ldots, T_5] 
\ / \ 
\bangle{T_1T_2^5 + T_3^3 + T_4^2}
$
\hspace{.5cm}
&
$(1,1,2,3,1)$
\hspace{.5cm}
&
$8$
\\
\midrule
2
&
$
\KK[T_1, \ldots, T_5] 
\ / \ 
\bangle{T_1T_2T_3^4 + T_4^3 + T_5^2}
$
&
$(1,1,1,2,3)$
&
$8$
\\
\midrule
3
&
$
\KK[T_1, \ldots, T_5] 
\ / \ 
\bangle{T_1T_2^2T_3^3 + T_4^3 + T_5^2}
$
&
$(1,1,1,2,3)$
&
$8$
\\
\midrule
4
&
$
\KK[T_1, \ldots, T_5] 
\ / \ 
\bangle{T_1T_2 + T_3T_4 + T_5^2}
$
&
$(1,1,1,1,1)$
&
$54$
\\
\midrule
5
&
$
\KK[T_1, \ldots, T_5] 
\ / \ 
\bangle{T_1T_2^2 + T_3T_4^2 + T_5^3}
$
&
$(1,1,1,1,1)$
&
$24$
\\
\midrule
6
&
$
\KK[T_1, \ldots, T_5] 
\ / \ 
\bangle{T_1T_2^3 + T_3T_4^3 + T_5^4}
$
&
$(1,1,1,1,1)$
&
$4$
\\
\midrule
7
&
$
\KK[T_1, \ldots, T_5] 
\ / \ 
\bangle{T_1T_2^3 + T_3T_4^3 + T_5^2}
$
&
$(1,1,1,1,2)$
&
$16$
\\
\midrule
8
&
$
\KK[T_1, \ldots, T_5] 
\ / \ 
\bangle{T_1T_2^5 + T_3T_4^5 + T_5^2}
$
&
$(1,1,1,1,3)$
&
$2$
\\
\midrule
9
&
$
\KK[T_1, \ldots, T_5] 
\ / \ 
\bangle{T_1T_2^5 + T_3^3T_4^3 + T_5^2}
$
&
$(1,1,1,1,3)$
&
$2$
\\
\bottomrule 
\end{longtable}
\end{center}
\end{theorem}


\subsection{$\KK^*$-surfaces}
Here we take a closer look at the 
first non-trivial examples of 
complexity one torus actions.
$\KK^*$-surfaces
are studied by many authors;
a classical reference is the work 
of Orlik and Wagreich~\cite{OrWa}.
Recall that a fixed point of a 
$\KK^*$-surface is said to be 
\begin{itemize}
\item
{\em elliptic\/} if it is isolated and lies in the 
closure of infinitely many $\KK^*$-orbits,
\item
{\em parabolic\/} if it belongs to a fixed point curve,
\item 
{\em hyperbolic\/} if it is isolated and lies in the 
closure of two $\KK^*$-orbits. 
\end{itemize}
These are in fact the only possible types of fixed points
for a normal $\KK^*$-surface. 
For normal projective $\KK^*$-surfaces $X$,
there is always a {\em source\/} $F^+ \subseteq X$ 
and a {\em sink\/} $F^- \subseteq X$.
They are characterized by the behaviour of
general points: there is a non-empty open set 
$U \subseteq X$ with 
$$ 
\lim_{t \to 0} t \mal x \ \in \ F^+,
\qquad
\lim_{t \to \infty} t \mal x \ \in \ F^-
\qquad
\text{ for all }
x \in U.
$$
The source can either consist of an elliptic fixed 
point or it is a curve of parabolic fixed points;
the same holds for the sink.
Any fixed point outside source or sink is hyperbolic.
Note that in Example~\ref{ex:blp1p1}, we have two curves 
of parabolic fixed points and three hyperbolic fixed points. 
To every smooth projective
$\KK^*$-surface $X$ having no elliptic fixed points
Orlik and Wagreich associated a graph of the 
following shape:
$$
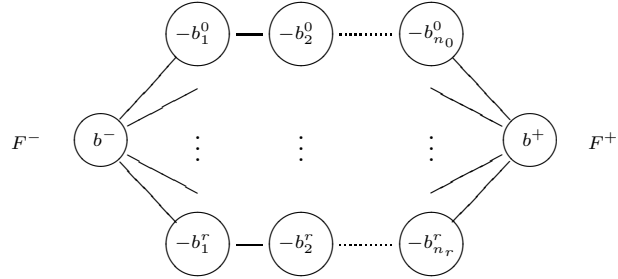
\entrymodifiers={++[o][F-]}
\def\objectstyle{\scriptstyle}
\xymatrix@-1.2pc{
*{}         
& *{}
& {-b^0_{1_{ \ }}} \ar@{-}[r]  
& {-b^0_{2_{ \ }}} \ar@{.}[rr] 
& *{}
& {-b^0_{n_0}}  
& *{}
& *{} 
\\
*{}
& *{}
& *{}
& *{}
& *{}
& *{}
& *{}
& *{}
\\
*{F^-}
&
{ \  b^-_{ \ }} 
\ar@{-}[uur] \ar@{-}[ur] \ar@{-}[ddr] \ar@{-}[dr] 
& *{\vdots}   
&*{\vdots}  
& *{}
& *{\vdots} 
& { \ b^+_{ \ }} 
\ar@{-}[ul] \ar@{-}[uul] \ar@{-}[ddl] \ar@{-}[dl] 
& *{F^+}
\\
*{}
& *{}
& *{}
& *{}
& *{}
& *{}
& *{}
& *{} 
\\
*{}
& *{}
& {-b_{1_{ \ }}^r} \ar@{-}[r]  
& {-b_{2_{ \ }}^r} \ar@{.}[rr] 
& *{}
& {-b_{n_r}^r} 
& *{}
& *{} 
\\
  }
$$
The two fixed point curves of $X$ occur as 
$F^-$ and $F^+$ in the graph.
The other vertices represent the invariant
irreducible contractible curves 
$D_{ij} \subseteq X$ different from $F^-$ and $F^+$.
The label $-b_{j}^i$ is the self intersection number
of $D_{ij}$, and two of the $D_{ij}$ are joined by an 
edge if and only if they have a common (fixed) point.
Every $D_{ij}$ is the closure of a non-trivial
$\KK^*$-orbit.

We show how to read off the Cox ring 
of a rational $X$ from its Orlik-Wagreich graph.
By~\cite[Sec.~3.5]{OrWa}, the order 
$l_{ij}$ of generic isotropy group 
of $D_{ij}$ equals the numerator of 
the canceled continued fraction
$$
b^i_1-\cfrac{1}{b^i_2-\cfrac{1}{\cdots -\cfrac{1}{b^i_{j-1}}}}
$$
Moreover, there is a canonical isomorphism 
$\PP_1 = Y = F^-$ identifying  $a_i \in Y$ 
with the point in $F^- \cap D_{i1}$.
In the terminology introduced before, we thus 
obtain the following.

\begin{theorem}
\label{OrWa2Cox}
Let $X$ be a smooth complete rational $\KK^*$-surface
without elliptic fixed points.
Then the assignments $S^{\pm} \mapsto 1_{F^{\pm}}$ and 
$T_{ij} \mapsto 1_{D_{ij}}$ define an isomorphism
\begin{eqnarray*}
\mathcal{R}(X)
&  \cong &
\KK[S^+,S^-,T_{ij}; \;  \; 0 \le i \le r, \; 1 \le j \le n_i] 
\ / \ 
\bangle{g_i; \; 0 \le i \le r-2}
\end{eqnarray*}
of $\Cl(X)$-graded rings, where the $\Cl(X)$-grading on the 
right hand side is defined by associating to $S^{\pm}$ the 
class of $F^{\pm}$ and to $T_{ij}$ the class of $D_{ij}$.
\end{theorem}

The Cox ring of a general rational  
$\KK^*$-surface $X$ is then computed as follows.
Blowing up the (possible) elliptic fixed points
suitably often gives a surface with two parabolic 
fixed point curves. Resolving the remaining 
singularities, we arrive at a smooth $\KK^*$-surface 
$\t{X}$ without parabolic fixed points, the so called 
{\em canonical resolution\/} of $X$; note that this 
is in general not the minimal resolution.
The Cox ring $\mathcal{R}(\t{X})$ is computed
as above.
For the Cox ring $\mathcal{R}(X)$, we need 
the divisors of the type $E_k$ and $D_{ij}$ in $X$ and 
the orders $l_{ij}$ of the generic isotropy groups
of the $D_{ij}$.
Each of these divisors is the image of 
a non-exceptional divisor of the same 
type in $\t{X}$; to see this for the $D_{ij}$,
note that $X_0$ is the open subset of $\t{X}_0$
obtained by removing the exceptional locus 
of $\t{X} \to X$ and thus $X_0/\KK^*$
is an open subset of $\t{X}_0/\KK^*$.
Moreover, by equivariance, 
the orders $l_{ij}$ in $X$ are the same 
as in $\t{X}$. 
Consequently, the Cox ring $\mathcal{R}(X)$
is obtained 
from $\mathcal{R}(\t{X})$ by removing
those generators that correspond to the exceptional 
curves arising from the resolution.

Now let us look at rational normal 
complete $\KK^*$-surfaces 
$X$ with the methods presented in the 
preceding subsection.
First recall that by finite generation 
of its Cox ring, $X$ is $\QQ$-factorial
and projective.
Moreover, $X$ arises from a ring $R(A,P)$ 
as in Construction~\ref{constr:RAP2Taction}.
As any surface, $X$ is the only 
representative of its small birational
class which means that no $\mathfrak{F}$-bunch
$\Phi$ needs to be specified in the defining 
data.

\begin{remark}
Let $(A,P)$ be data as in~\ref{constr:RAP}.
In order that $R(A,P)$ defines a surface~$X$,
we necessarily have $s=1$ and for $m$ we have 
the following three possibilities
\begin{itemize}
\item 
$m = 0$. This holds if and only if $X$ has two
elliptic fixed points. In this case $d'$ is empty.
\item
$m = 1$. This holds if and only if $X$ has one
fixed point curve and one 
elliptic fixed point. In this case, we may 
assume $d' = (-1)$.
\item
$m = 2$. This holds if and only if $X$ has two
elliptic fixed points. In this case, we may 
assume $d' = (1,-1)$.
\end{itemize}
In the fan $\Sigma$ of the canonical ambient 
toric variety $Z$ of $X$, the 
parabolic fixed point curves correspond to 
rays through $(0,\ldots,0,\pm 1)$,
the elliptic fixed points to
full dimensional cones $\sigma^{\pm}$ and 
the hyperbolic fixed points to two-dimensional
cones. 
\end{remark}

We demonstrate the concrete working with
$\KK^*$-surfaces by resolving singularities 
in two examples;
for a detailed general treatment we refer 
to~\cite{He}.
Similar to the canonical resolution
of Orlik and Wagreich, the first step
is resolving singular elliptic 
fixed points by inserting rays through 
$(0,\ldots,0,\pm 1)$.
The remaining singularities are then 
resolved by regular subdivision of 
two-dimensional singular cones of 
$\Sigma$.
The first example is our surface~\ref{ex:delpezzo}. 

\begin{example}
In the setting of~\ref{constr:RAP}, let
$r = 2$, set $n_0=2$, $n_1=1$, $n_2=2$ 
and $m=0$. For $A$ choose the vectors 
$(-1,0)$, $(1,-1)$ and $(0,1)$.
Consider the matrix
\begin{eqnarray*}
P 
& = &
\left(
\begin{array}{rrrrr}
-1 & -1 & 2 & 0 & 0
\\
-1 & -1 & 0 & 1 & 1 
\\
-1 & 0 & 1 & -1 & 0
\end{array}
\right).
\end{eqnarray*}
Then the defining equation of $\b{X} \subseteq \KK^5$ 
is $g = T_{01}T_{02} + T_{11}^2 + T_{21}T_{22}$.
Observe that the Gale dual matrix is
\begin{eqnarray*}
Q 
& = &
\left(
\begin{array}{rrrrr}
1 & -1 & 0 & -1 & 1
\\
1 & 1 & 1 & 0 & 2 
\end{array}
\right).
\end{eqnarray*}
Thus, $X$ is the del Pezzo surface we 
already encountered in~\ref{ex:delpezzo}.
Consider the canonical toric ambient variety $Z$.
Its fan $\Sigma$ looks as follows:
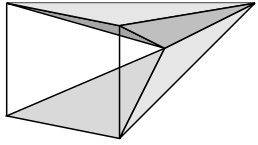
\begin{center}
\begin{picture}(0,0)%
\includegraphics{A2delpambientfan.pstex}%
\end{picture}%
\setlength{\unitlength}{1243sp}%
\begingroup\makeatletter\ifx\SetFigFontNFSS\undefined%
\gdef\SetFigFontNFSS#1#2#3#4#5{%
  \reset@font\fontsize{#1}{#2pt}%
  \fontfamily{#3}\fontseries{#4}\fontshape{#5}%
  \selectfont}%
\fi\endgroup%
\begin{picture}(5016,2766)(3118,-4144)
\end{picture}%

\end{center}
More precisely, in terms of the columns 
$v_{01}$, $v_{02}$, $v_{11}$, $v_{21}$, $v_{22}$ 
of $P$, the maximal cones of $\Sigma$ are 
$$ 
\sigma^- \ := \ \cone(v_{01},v_{11},v_{21}),
\qquad
\sigma^+ \ := \ \cone(v_{02},v_{12},v_{22}),
$$
$$
\tau_{012} \ := \ \cone(v_{01},v_{02}),
\qquad
\tau_{212} \ := \ \cone(v_{21},v_{22}).
$$
As we already observed, $X$ comes with 
one singularity $x_0 \in X$;
it is the toric fixed point corresponding 
to the cone $\sigma^-$. Set
$$ 
v_1 \ := \ (0,0,-1),
\qquad
v_{10} \ := \ (1,0,0).
$$
Inserting these vectors as columns 
at the right places of the 
matrix $P$ gives the describing matrix 
for the resolution $\t{X}$ of $X$:
\begin{eqnarray*}
\t{P }
& = &
\left(
\begin{array}{rrrrrrr}
-1 & -1 & 1 & 2 & 0 & 0 & 0
\\
-1 & -1 & 0 & 0 & 1 & 1 & 0
\\
-1 & 0 & 0 & 1 & -1 & 0 & -1
\end{array}
\right).
\end{eqnarray*}
This allows us to read off the defining 
relation of the Cox ring; 
the degrees of the generators are the 
columns of the Gale dual matrix $\t{Q}$.
Concretely, we obtain:
\begin{eqnarray*}
\mathcal{R}(X)
& = & 
\KK[T_{ij},S; \; i=0,1,2 \; j = 1,2] 
\ / \ 
\bangle{T_{01}T_{02} + T_{11}T_{12}^2 + T_{21}T_{22}},
\\[2ex]
\t{Q}
& := & 
\left(
\begin{array}{rrrrrrr}
0 & 1 & -1 & 1 & 1 & 0 & 0
\\
1 & 0 & -1 & 1 & 0 & 1 & 0
\\
0 & -1 & -1 & 0 & 0 & -1 & 0
\\
0 & 0 & 0 & 0 & -1 & 1 & 1
\end{array}
\right).
\end{eqnarray*}
The self intersection numbers of the 
curves
corresponding to $v_1$ and $v_{10}$
are computed as in~\ref{rem:compint} 
and both equal $-2$. 
Thus, the singularity of $X$ is of type~$A_2$.
\end{example}

In the second example, we compute the
Cox ring of the minimal resolution of 
the $E_6$ cubic surface from~\ref{ex:E6cubic};
the result was first obtained by Hassett 
and Tschinkel, without using the 
$\KK^*$-action~\cite[Theorem~3.8]{HaTs}.

\begin{example}
\label{ex:E6cubicresol}
In the setting of~\ref{constr:RAP}, let
$r = 2$, set $n_0=2$, $n_1=1$, $n_2=1$ 
and $m=0$. For $A$ choose the vectors 
$(-1,0)$, $(1,-1)$ and $(0,1)$.
Consider the matrix
\begin{eqnarray*}
P 
& = &
\left(
\begin{array}{rrrr}
-1 & -3 & 3 & 0 
\\
-1 & -3 & 0 & 2  
\\
-1 & -2 & 1 & 1 
\end{array}
\right).
\end{eqnarray*}
The resulting surface $X$ is the $E_6$
cubic surface~\ref{ex:E6cubic}.
Its Cox ring and degree matrix are given 
as
\begin{eqnarray*}
\mathcal{R}(X)
& = & 
\KK[T_{01},T_{02},T_{11},T_{21}]
\ / \  
\bangle{T_{01}T_{02}^3 + T_{11}^3 + T_{21}^2}
\\[2ex]
Q 
& := & 
\left(
\begin{array}{rrrr}
3 & 1 & 2 & 3
\end{array}
\right).
\end{eqnarray*}
Let us look at the fan of the toric ambient variety.
In terms of the columns $v_{01},v_{02},v_{11},v_{21}$
of $P$, its maximal cones are
$$ 
\sigma^- \ := \ \cone(v_{01},v_{11},v_{21}),
\qquad
\sigma^+ \ := \ \cone(v_{02},v_{11},v_{22}),
\qquad
\tau_{012} \ := \ \cone(v_{01},v_{02}).
$$
The toric fixed point corresponding to 
$\sigma^+$ is the singularity $x_0 \in X$.
It is singular for two reasons: firstly 
$\sigma^+$ is singular, secondly the total 
coordinate space
\begin{eqnarray*}
\b{X} & = & V(T_{01}T_{02}^3 + T_{11}^3 + T_{21}^2)
\end{eqnarray*}
is singular along the fiber 
$\KK^* \times 0 \times 0 \times 0$
over $x_0$.
Subdividing along $(0,0,1)$ and further resolving 
gives
\begin{eqnarray*}
\t{P }
& = &
\left(
\begin{array}{rrrrrrrrrr}
-1 & -3 & -2 & -1 & 3 & 2 & 1 & 0 & 0 & 0 
\\
-1 & -3 & -2 & -1 & 0 & 0 & 0 & 2 & 1 & 0 
\\
-1 &  -2 & -1 & 0 & 1 & 1 & 1 & 1 & 1 & 1
\end{array}
\right).
\end{eqnarray*}
Note that inserting first $(0,0,1)$ is part of 
the first step of the canonical resolution.
The Cox ring of the resolution and its degree 
matrix are thus given by
\begin{eqnarray*}
\mathcal{R}(\t{X})
& = & 
\frac{
\KK[T_{01},T_{02},T_{03},T_{04},T_{11},T_{12},T_{13},T_{21},T_{21},S]
}
{
\bangle{T_{01}T_{02}^3T_{03}^2T_{04} 
        + T_{11}^3T_{12}^2T_{13}
        + T_{21}^2T_{22}}
},
\\[2ex]
\t{Q}
& = & 
\left(
\begin{array}{rrrrrrrrrr}
-1 & 1 & -1 & 0 & 0 & 0 & 0 & 0 & 0 & 0
\\
1 & 0 & -1 & 1 & 0 & 0 & 0 & 0 & 0 & 0
\\
0 & 1 & 0 & -1 & 0 & 1 & 0 & 1 & 0 & 0
\\
0 & 0 & -1 & 0 & -1 & 0 & 1 & -1 & 0 & 0
\\
1 & 0 & 0 & 0 & 0 & 1 & -1 & 0 & 1 & 0
\\
-1 & 0 & 0 & 0 & -1 & 1 & 0 & 0 & -1 & 0
\\
1 & 0 & 0 & -1 & 0 & 0 & 0 & 0 & 0 & 1
\end{array}
\right).
\end{eqnarray*}
\end{example}

\begin{remark}
Inserting the rays through $(0,\ldots,0,\pm1)$ as
(partially) done in the preceding two examples 
is a ``tropicalization step'': the two rays arise
when we intersect the fan of the toric ambient 
variety with the tropicalization ${\rm trop}(X)$
which, in the example case, is the union of the 
one codimensional cones of the normal fan of 
the Newton polytope of the defining equation.
\end{remark}

Combining the description in terms of 
data $(A,P)$ with the language of bunched 
rings makes $\KK^*$-surfaces an easily 
accessible class of examples. 
This allows among other things explicit 
classifications. 
For example, the Gorenstein del Pezzo 
$\KK^*$-surfaces are classified by these
methods in~\cite{He}; part of the results
has been obtained by other methods 
in~\cite{Der} and~\cite{Suss,HaSu}.

\begin{theorem}
The following table lists 
Cox ring $\mathcal{R}(X)$
and the singularity type $\mathrm S(X)$ 
of the non-toric Gorenstein Fano 
$\KK^*$-surfaces $X$ of Picard number one.
\begin{center}
\begin{longtable}[htbp]{cccc}
\toprule
$\mathcal{R}(X)$ 
&
$\Cl(X)$
& 
$(w_1,\ldots, w_r)$
& 
$\mathrm S(X)$
\\
\midrule
$
\KK[T_1, T_2 , T_3,S_1] 
\ / \ 
\bangle{T_1^2 + T_2^2 + T_3^2}
$
&
$\ZZ\oplus \ZZ / 2\ZZ\oplus \ZZ / 2\ZZ$
&
$
\left(
\begin{smallmatrix}
1 & 1 & 1 & 1\\
\bar 1 & \bar 0 & \bar 1 & \bar 0\\
\bar 1 & \bar 1 & \bar 0 & \bar 0
\end{smallmatrix}
\right)
$
&
$D_43A_1$
\\
\cmidrule{1-4}
$
\KK[T_1, \ldots, T_4] 
\ / \ 
\bangle{T_1T_2 + T_3^2 + T_4^2}
$
&
$\ZZ\oplus \ZZ / 4\ZZ$
&
$
\left(
\begin{smallmatrix}
1 & 1 & 1 & 1\\
\bar 1 & \bar 3 & \bar 2  & \bar 0
\end{smallmatrix}
\right)
$
&
$2A_3A_1$
\\
\cmidrule{1-4}
$
\KK[T_1, \ldots, T_4] 
\ / \ 
\bangle{T_1T_2 + T_3^2 + T_4^3}
$
&
$\ZZ$
&
$
\left(
\begin{smallmatrix}
1 & 5 & 3 & 2
\end{smallmatrix}
\right)
$
&
$A_4$
\\
\cmidrule{1-4}
$
\KK[T_1, \ldots, T_4] 
\ / \ 
\bangle{T_1T_2 + T_3^2 + T_4^4}
$
&
$\ZZ \oplus \ZZ / 2\ZZ$
&
$\left(
\begin{smallmatrix}
1 & 3 & 2 & 1\\
\bar 1 & \bar 1 & \bar 1 & \bar 0
\end{smallmatrix}
\right)
$
&
$A_5A_1$
\\
\cmidrule{1-4}
$
\KK[T_1, \ldots, T_4] 
\ / \ 
\bangle{T_1T_2 + T_3^3 + T_4^3}
$
&
$\ZZ\oplus \ZZ / 3\ZZ$
&
$
\left
(\begin{smallmatrix}
1 & 2 & 1 & 1\\
\bar 1 & \bar 2 &  \bar 2 & \bar 0
\end{smallmatrix}
\right)
$
&
$A_5A_2$
\\
\cmidrule{1-4}
$
\KK[T_1, \ldots, T_4] 
\ / \ 
\bangle{T_1T_2^2 + T_3^3 + T_4^2}
$
&
$\ZZ$
&
$
\left(
\begin{smallmatrix}
4 & 1 & 2 & 3
\end{smallmatrix}
\right)
$
&
$D_5$
\\
\cmidrule{1-4}
$
\KK[T_1, \ldots, T_4] 
\ / \ 
\bangle{T_1^2T_2 + T_3^2 + T_4^4}
$
&
$\ZZ\oplus \ZZ / 2\ZZ$
&
$\left(
\begin{smallmatrix}
1 & 2 & 2 & 1\\
\bar 1 & \bar 0 & \bar 1  & \bar 0
\end{smallmatrix}
\right)
$
&
$D_6A_1$
\\
\cmidrule{1-4}
$
\KK[T_1, \ldots, T_4] 
\ / \ 
\bangle{T_1T_2^2 + T_3^3 + T_4^3}
$
&
$\ZZ\oplus \ZZ / 3\ZZ$
&
$
\left(
\begin{smallmatrix}
1 & 1 & 1 & 1\\
\bar 1 & \bar 1 & \bar 2 & \bar 0
\end{smallmatrix}
\right)
$
&
$E_6A_2$
\\
\cmidrule{1-4}
$
\KK[T_1, \ldots, T_4] 
\ / \ 
\bangle{T_1T_2^3 + T_3^3 + T_4^2}
$
&
$\ZZ$
&
$
\left(
\begin{smallmatrix}
3 & 1 & 2 & 3
\end{smallmatrix}
\right)
$
&
$E_6$
\\
\cmidrule{1-4}
$
\KK[T_1, \ldots, T_4] 
\ / \ 
\bangle{T_1T_2^3 + T_3^4 + T_4^2}
$
&
$\ZZ\oplus \ZZ / 2\ZZ$
&
$
\left(
\begin{smallmatrix}
1 & 1 & 1 & 2\\
\bar 0 & \bar 0 & \bar 1 & \bar 1
\end{smallmatrix}
\right)
$
&
$E_7A_1$
\\
\cmidrule{1-4}
$
\KK[T_1, \ldots, T_4] 
\ / \ 
\bangle{T_1T_2^4 + T_3^3 + T_4^2}
$
&
$\ZZ$
&
$
\left(
\begin{smallmatrix}
2 & 1 & 2 & 3
\end{smallmatrix}
\right)
$
&
$E_7$
\\
\cmidrule{1-4}
$
\KK[T_1, \ldots, T_4] 
\ / \ 
\bangle{T_1T_2^5 + T_3^3 + T_4^2}
$
&
$\ZZ$
&
$
\left(
\begin{smallmatrix}
1 & 1 & 2 & 3
\end{smallmatrix}
\right)
$
&
$E_8$
\\
\cmidrule{1-4}
$
\KK[T_1, \ldots, T_5] 
\ / \ 
\bangle{
\begin{smallmatrix}
T_1T_2 + T_3^2 + T_4^2,\\ \lambda T_3^2+T_4^2+T_5^2
\end{smallmatrix}}
$
&
$\ZZ\oplus \ZZ / 2\ZZ\oplus \ZZ / 2\ZZ$
&
$\left(
\begin{smallmatrix}
1 & 1 & 1 & 1 & 1\\
\bar 1 &\bar 1 & \bar 0 & \bar 1 &  \bar 0\\
\bar 0 & \bar 0 & \bar 1 & \bar 1 &  \bar 0
\end{smallmatrix}
\right)$
&
$2D_4$
\\
\bottomrule 
\end{longtable}
\end{center}
\end{theorem}

\begin{theorem}
The following table lists 
Cox ring $\mathcal{R}(X)$
and the singularity type $\mathrm S(X)$ 
of the non-toric Gorenstein Fano 
$\KK^*$-surfaces $X$ of Picard number two.
\begin{center}
\begin{longtable}[htbp]{cccc}
\toprule
$\mathcal{R}(X)$ 
& 
$\Cl(X)$
&
$(w_1,\ldots, w_r)$
& 
$\mathrm S(X)$
\\
\midrule
$
\KK[T_1, \ldots,T_4,S_1] 
\ / \ 
\bangle{T_1T_2 + T_3^2 + T_4^2}
$
&
$\ZZ^2 \oplus \ZZ/2\ZZ$
&
$
\left(
\begin{smallmatrix}
1 & 1 & 1 & 1 & 1\\
1 & -1 & 0 & 0 & 1\\
\bar 1 & \bar 1 & \bar 1 & \bar 0  & \bar 0
\end{smallmatrix}
\right)$
&
$A_32A_1$
\\
\cmidrule{1-4}
$
\KK[T_1, \ldots,T_5] 
\ / \ 
\bangle{T_1T_2 + T_3T_4 + T_5^2}
$
&
$\ZZ^2$
&
$
\left(
\begin{smallmatrix}
1 & 1 & 1 & 1 & 1\\
-1 & 1 & 2 & -2 & 0
\end{smallmatrix}
\right)
$
&
$2A_2A_1$
\\
\cmidrule{1-4}
$
\KK[T_1, \ldots, T_5] 
\ / \ 
\bangle{T_1T_2 + T_3T_4 + T_5^2}
$
&
$\ZZ^2$
&
$
\left(
\begin{smallmatrix}
1 & 3 & 1 & 3 & 2\\
1 & 1 & 0 & 2 & 1
\end{smallmatrix}
\right)
$
&
$A_2$
\\
\cmidrule{1-4}
$
\KK[T_1, \ldots, T_5] 
\ / \ 
\bangle{T_1T_2 + T_3T_4 + T_5^3}
$
&
$\ZZ^2$
&
$
\left(
\begin{smallmatrix}
1 & 2 & 1 & 2 & 1\\
1 & -1 & -1 & 1 & 0
\end{smallmatrix}
\right)
$
&
$A_1A_3$
\\
\cmidrule{1-4}
$
\KK[T_1, \ldots, T_5] 
\ / \ 
\bangle{T_1T_2 + T_3^2T_4 + T_4^2}
$
&
$\ZZ^2$
&
$
\left(
\begin{smallmatrix}
1 & 3 & 1 & 2 & 2\\
0 & 2 & 1 & 0 & 1
\end{smallmatrix}
\right)
$
&
$A_3$
\\
\cmidrule{1-4}
$
\KK[T_1, \ldots, T_5] 
\ / \ 
\bangle{T_1T_2 + T_3^2T_4 + T_5^3}
$
&
$\ZZ^2$
&
$
\left(
\begin{smallmatrix}
1 & 2 & 1 & 1 & 1\\
-1 & 1 & 1 & -2 & 0
\end{smallmatrix}
\right)
$
&
$A_4A_1$
\\
\cmidrule{1-4}
$
\KK[T_1, \ldots, T_5] 
\ / \ 
\bangle{T_1T_2 + T_3^3T_4 + T_5^2}
$
&
$\ZZ^2$
&
$
\left(
\begin{smallmatrix}
1 & 3 & 1 & 1 & 2\\
-1 & -1 & 0 & -2 & -1
\end{smallmatrix}
\right)
$
&
$A_4$
\\
\cmidrule{1-4}
$
\KK[T_1, \ldots, T_5] 
\ / \ 
\bangle{T_1T_2^2 + T_3T_4^2 + T_5^2}
$
&
$\ZZ^2$
&
$
\left(
\begin{smallmatrix}
2 & 1 & 2 & 1 & 2\\
0 & -1 & -2 & 0 & -1
\end{smallmatrix}
\right)
$
&
$D_4$
\\
\cmidrule{1-4}
$
\KK[T_1, \ldots, T_5] 
\ / \ 
\bangle{T_1T_2^2 + T_3T_4^2 + T_5^3}
$
&
$\ZZ^2$
&
$
\left(
\begin{smallmatrix}
1 & 1 & 1 & 1 & 1\\
2 & -1 & -2 & 1 & 0
\end{smallmatrix}
\right)
$
&
$D_5A_1$
\\
\cmidrule{1-4}
$
\KK[T_1, \ldots, T_5] 
\ / \ 
\bangle{T_1T_2^3 + T_3T_4^2 + T_5^2}
$
&
$\ZZ^2$
&
$
\left(
\begin{smallmatrix}
1 & 1 & 2 & 1 & 2\\
1 & -1 & -2 & 0 & -1
\end{smallmatrix}
\right)
$
&
$D_5$
\\
\cmidrule{1-4}
$
\KK[T_1, \ldots, T_5] 
\ / \ 
\bangle{T_1T_2^3 + T_3T_4^3 + T_5^2}
$
&
$\ZZ^2$
&
$
\left(
\begin{smallmatrix}
1 & 1 & 1 & 1 & 2\\
1 & -1 & -2 & 0 & -1
\end{smallmatrix}
\right)
$
&
$E_6$
\\
\cmidrule{1-4}
$
\KK[T_1, \ldots, T_6] 
\ / \ 
\bangle{
\begin{smallmatrix}
T_1T_2 + T_3T_4 + T_5^2,\\ \lambda T_3T_4+T_5^2+T_6^2
\end{smallmatrix}}
$
&
$\ZZ^2\oplus \ZZ/2\ZZ$
&
$
\left(
\begin{smallmatrix}
1 & 1 & 1 & 1 & 1 & 1\\
-1 & 1 & 1 & -1 & 0 & 0\\
\bar 0 & \bar 0 & \bar 1 & \bar 1 & \bar 1 &  \bar 0
\end{smallmatrix}
\right)
$
&
$2A_3$
\\
\bottomrule 
\end{longtable}
\end{center}
\end{theorem}

\begin{theorem}
The following table lists 
Cox ring $\mathcal{R}(X)$ 
and the singularity type $\mathrm S(X)$ 
of the non-toric Gorenstein Fano 
$\KK^*$-surfaces $X$ of Picard number three.
\begin{center}
\begin{longtable}[htbp]{cccc}
\toprule
$\mathcal{R}(X)$ 
&
$\Cl(X)$
& 
$(w_1,\ldots, w_r)$
& 
$\mathrm S(X)$
\\
\midrule
$
\KK[T_1, \ldots,T_5,S_1] 
\ / \ 
\bangle{T_1T_2 + T_3T_4 + T_5^2}
$
&
$\ZZ^3$
&
$\left(
\begin{smallmatrix}
1 & 1 & 1 & 1 & 1 & 1 \\
1 & -1 & 0 & 0 & 0 & 1\\
1 & -1 & -1 & 1 & 0 & 0
\end{smallmatrix}
\right)
$
&
$A_1A_2$
\\
\cmidrule{1-4}
$
\KK[T_1, \ldots,T_6] 
\ / \ 
\bangle{T_1T_2 + T_3T_4 + T_5T_6}
$
&
$\ZZ^3$
&
$
\left(
\begin{smallmatrix}
1 & 1 & 1 & 1 & 1 & 1 \\
1 & 0 & 0 & 1 & 0 & 1\\
0 & 0 & 1 & -1 & -1 & 1
\end{smallmatrix}
\right)
$
&
$3A_1$
\\
\cmidrule{1-4}
$
\KK[T_1, \ldots,T_6] 
\ / \ 
\bangle{T_1T_2 + T_3T_4 + T_5T_6}
$
&
$\ZZ^3$
&
$
\left(
\begin{smallmatrix}
1 & 2 & 1 & 2 & 1 & 2 \\
0 & 1 & 0 & 1 & 1 & 0\\
1 & -1 & -1 & 1 & 0 & 0 
\end{smallmatrix}
\right)
$
&
$A_1$
\\
\cmidrule{1-4}
$
\KK[T_1, \ldots, T_6] 
\ / \ 
\bangle{T_1T_2 + T_3T_4 + T_5^2T_6}
$
&
$\ZZ^3$
&
$
\left(
\begin{smallmatrix}
1 & 2 & 1 & 2 & 1 & 1 \\
1 & 0 & 0 & 1 & 0 & 1\\
1 & -1 & -1 & 1 & 0 & 0 
\end{smallmatrix}
\right)
$
&
$A_2$
\\
\cmidrule{1-4}
$
\KK[T_1, \ldots, T_6] 
\ / \ 
\bangle{T_1T_2 + T_3T_4^2 + T_5T_6^2}
$
&
$\ZZ^3$
&
$
\left(
\begin{smallmatrix}
2 & 1 & 1 & 1 & 1 & 1 \\
0 & 1 & 1 & 0 & 1 & 0\\
1 & 0 & -1 & 1 & 1 & 0 
\end{smallmatrix}
\right)
$
&
$A_3$
\\
\cmidrule{1-4}
$
\KK[T_1, \ldots, T_6] 
\ / \ 
\bangle{T_1T_2^2 + T_3T_4^2 + T_5T_6^2}
$
&
$\ZZ^3$
&
$
\left(
\begin{smallmatrix}
1 & 1 & 1 & 1 & 1 & 1 \\
1 & -1 & -1 & 0 & -1 & 0\\
1 & 0 & -1 & 1 & 1 & 0 
\end{smallmatrix}
\right)
$
&
$D_4$
\\
\cmidrule{1-4}
$
\KK[T_1, \ldots, T_7] 
\ / \ 
\bangle{
\begin{smallmatrix}
T_1T_2 + T_3T_4 + T_5T_6,\\ \lambda T_3T_4+T_5T_6+T_7^2
\end{smallmatrix}}
$
&
$\ZZ^3$
&
$
\left(
\begin{smallmatrix}
1 & 1 & 1 & 1 & 1 & 1 & 1 \\
0 & 0 & 1 & -1 & -1 & 1 & 0\\
-1 & 1 & 1 & -1 & 0 & 0 & 0
\end{smallmatrix}
\right)
$
&
$2A_2$
\\
\bottomrule 
\end{longtable}
\end{center}
\end{theorem}

\begin{theorem}
There is one non-toric Gorenstein Fano 
$\KK^*$-surface of Picard number four;
its Cox ring $\mathcal{R}(X)$ 
and singularity type $\mathrm S(X)$ 
is given below.
\begin{center}
\begin{longtable}[htbp]{cccc}
\toprule
$\mathcal{R}(X)$ 
&
$\Cl(X)$
& 
$(w_1,\ldots, w_r)$
& 
$\mathrm S(X)$
\\
\midrule
$
\KK[T_1, \ldots,T_6,S_1] 
\ / \ 
\bangle{T_1T_2 + T_3T_4 + T_5T_6}
$
&
$\ZZ^4$
&
$
\left(
\begin{smallmatrix}
1 & 1 & 1 & 1 & 1 & 1 & 1\\
1 & 0 & 0 & 1 & 0 & 1 & 0\\
0 & 1 & 0 & 1 & 1 & 0 & 0\\
1 & -1 & -1 & 1 & 0 & 0 & 0
\end{smallmatrix}
\right)
$
&
$A_1$
\\
\bottomrule 
\end{longtable}
\end{center}
\end{theorem}


\end{document}